\documentclass[a4paper,12pt]{article}
\usepackage{amsfonts}
\usepackage{amssymb,epic,eepic}
\usepackage{epsfig}
\begin{document}

\title{Nonlinear hyperbolic equations \\ in surface theory: \\
integrable discretizations \\and approximation results}
\author{A.I.\,Bobenko
 \thanks{ E--Mail: {\tt bobenko}@{\tt math.tu-berlin.de}}
\and D.\,Matthes
 \thanks{ E--Mail: {\tt matthes}@{\tt math.tu-berlin.de}}
\and Yu.B.\,Suris
 \thanks{ E--Mail: {\tt suris}@{\tt sfb288.math.tu-berlin.de}}}
\date{Institut f\"ur Mathematik, Technische Universit\"at Berlin,\\
 Str. des 17. Juni 136, 10623 Berlin, Germany.}

\newcommand{\ds}{\displaystyle}
\newcommand{\be}{\begin{equation}}
\newcommand{\ee}{\end{equation}}
\newcommand{\beq}{\begin{eqnarray*}}
\newcommand{\eeq}{\end{eqnarray*}}
\newcommand{\beqn}{\begin{eqnarray}}
\newcommand{\eeqn}{\end{eqnarray}}
\newcommand{\eps}{\epsilon}
\newcommand{\gz}{{\mathbb Z}}
\newcommand{\rz}{{\mathbb R}}
\newcommand{\cz}{{\mathbb C}}
\newcommand{\bfe}{{\mathbf e}}
\newcommand{\bfEps}{\mbox{\boldmath$\eps$}}
\newcommand{\bfeps}{\mbox{\small\boldmath$\eps$}}
\newcommand{\cD}{{\cal D}}
\newcommand{\cE}{{\cal E}}
\newcommand{\cF}{{\cal F}}
\newcommand{\cG}{{\cal G}}
\newcommand{\cK}{{\cal K}}
\newcommand{\cL}{{\cal L}}
\newcommand{\cO}{{\cal O}}
\newcommand{\cP}{{\cal P}}
\newcommand{\cU}{{\cal U}}
\newcommand{\cV}{{\cal V}}
\newcommand{\cW}{{\cal W}}
\newcommand{\cX}{{\cal X}}
\newcommand{\cY}{{\cal Y}}

\newcommand{\vfree}{\vspace{0pt}\\}
\newcommand{\nn}{\Vert}
\newcommand{\OH}{{\cal O}}

\newtheorem{thm}{Theorem}
\newtheorem{dfn}{Definition}
\newtheorem{lem}{Lemma}
\newtheorem{prp}{Proposition}
\newtheorem{cor}{Corollary}
\def\itbf{\itshape\bfseries}

\def \BBox{\hspace{1mm}\vrule height6pt width5.5pt depth0pt \hspace{6pt}}

\maketitle

\begin{abstract} A numerical scheme is developed for solution of
the Goursat problem for a class of nonlinear hyperbolic systems with
an arbitrary number of independent variables. Convergence results are proved
for this difference scheme. These results are applied to hyperbolic systems
of differential--geometric origin, like the sine--Gordon equation
describing the surfaces of the constant negative Gaussian curvature
($K$--surfaces). In particular, we prove the convergence of discrete
$K$--surfaces and their B\"acklund transformations to their continuous
counterparts. This puts on a firm basis the generally accepted belief (which
however remained unproved untill this work) that the classical differential
geometry of integrable classes of surfaces and the classical theory of
transformations of such surfaces may be obtained from a unifying
multi--dimensional discrete theory by a refinement of the coordinate
mesh--size in some of the directions.\end{abstract}

\newpage
\section{Introduction}

The development of the classical differential geometry led to
introduction and studying various classes of surfaces which are of
interest both for the internal differential--geometric reasons and
for application in other sciences. We mention here minimal
surfaces, constant curvature surfaces, isothermic surfaces, to
name just a few general classes. A rich theory of such surface
classes is, to a large extent, a classical heritage; more
recently, also numerical methods appeared, stimulated by
applications in sciences and in scientific computation, including
visualization. Mostly, these numerical methods are of a
variational nature, applicable to problems described by elliptic
partial differential equations, like the Plateau problem in the
theory of minimal surfaces (see, e.g., \cite{PP, Hi, DH}). It
seems that no general numerical methods have been developed for
objects of the differential geometry described by hyperbolic
differential equations, like the surfaces with constant negative
Gaussian curvature.

\begin{figure}[ht]
  \begin{center}
  \begin{picture}(0,0)%
\includegraphics{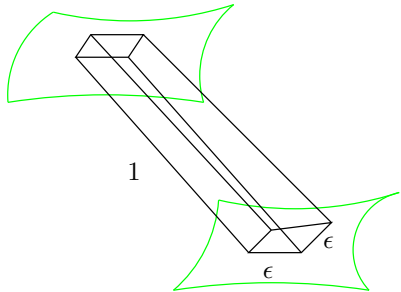}%
\end{picture}%
\setlength{\unitlength}{4144sp}%
\begingroup\makeatletter\ifx\SetFigFont\undefined%
\gdef\SetFigFont#1#2#3#4#5{%
  \reset@font\fontsize{#1}{#2pt}%
  \fontfamily{#3}\fontseries{#4}\fontshape{#5}%
  \selectfont}%
\fi\endgroup%
\begin{picture}(2359,1726)(260,-1464)
\put(991,-781){\makebox(0,0)[lb]{\smash{\SetFigFont{10}{12.0}{\rmdefault}{\mddefault}{\updefault}
$1$
}}}
\put(2161,-1186){\makebox(0,0)[lb]{\smash{\SetFigFont{10}{12.0}{\rmdefault}{\mddefault}{\updefault}
$\epsilon$
}}}
\put(1801,-1366){\makebox(0,0)[lb]{\smash{\SetFigFont{10}{12.0}{\rmdefault}{\mddefault}{\updefault}
$\epsilon$
}}}
\end{picture}
  \caption{Surfaces and their transformations as a limit
  of multidimensional lattices}  \label{Limit}
  \end{center}
\end{figure}

The characteristic property of various special classes of surfaces studied
by the classical differential geometry turns out to be their integrability,
in the sense of the modern theory of solitons. One of the manifestations of
integrability is the existence of a rich transformations theory,
sometimes unified under the names of the Darboux--B\"acklund
transformations. Classically, the theory of surfaces and the theory
of their transformations were dealt with separately
to a large extent. Recently, it became clear that both theories can
be unified in the framework of the discrete differential geometry (cf.
\cite{S, BP2}). In this framework, multidimensional lattices with certain
geometrical properties become the basic mathematical structures.
As the lattice becomes more and more dense in some of the
coordinates directions (the mesh size $\epsilon\to 0$), it
approximates the smooth surface. The directions, where the mesh
size remains constant, correspond to the transformations of smooth
surfaces (see Fig. \ref{Limit}).

\begin{center}
\begin{figure}[htbp]
\centerline{\epsffile{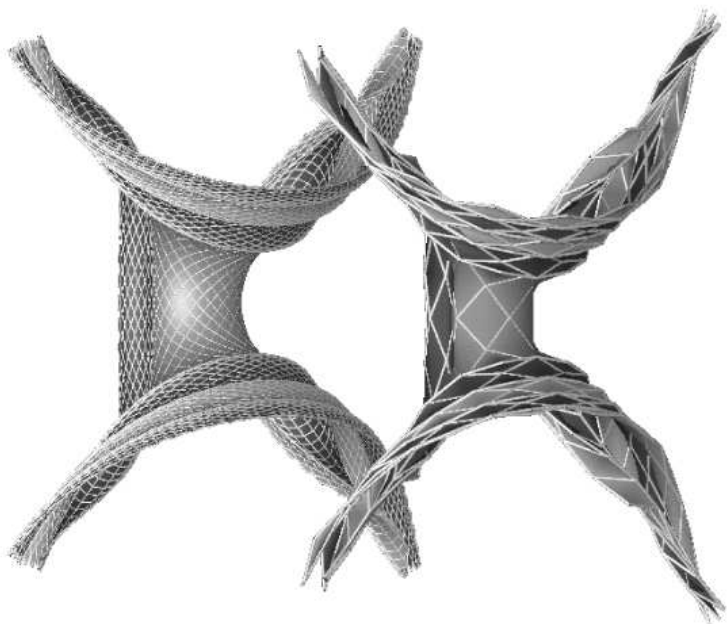}}
\caption{A countinous and a discrete Amsler surfaces} 
\label{Amsler}
\end{figure}
\end{center}

The discrete differential geometry is nowadays a flourishing area
which parallels to a large extent its classical (continuous)
counterpart. Many important classes of surfaces have been
discretized up to now, see a review in \cite{BP2}. Their
properties are well understood. Fig.\ref{Amsler} \footnote{We are 
    thankful to Tim Hoffmann for producing these figures.} shows an example
of a continuous Amsler surface (which is a surface with constant
negative Gaussian curvature) and its discrete analog. The
characteristic property of discrete surfaces
$F:(\eps\gz)^2\mapsto\rz^3$ with constant negative Gaussian
curvature is that for each $(x,y)\in(\eps\gz)^2$ the five points
$F(x,y)$ and $F(x\pm\eps,y\pm\eps)$ lie in a plane. The subclass
of Amsler surfaces is singled out by the condition that the
surface should contain two straight lines.

\epsfxsize=300pt

\begin{center}
\begin{figure}[htbp]
\centerline{\epsffile{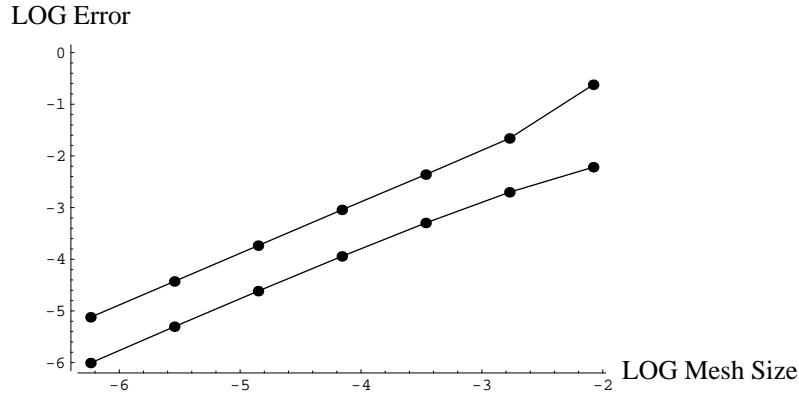}}
\caption{Maximum norm of the error $F^\eps -F$ vs. mesh size $\eps$.
$F:[0,r]\times[0,r]\rightarrow{\mathbb R}^3$ is a smooth K-surface,
$F^\eps:[0,r]^{\eps}\times[0,r]^\eps\rightarrow{\mathbb R}^3$ are
discrete K-surfaces. (See section \ref{Sect 2d} for notations.) The
lower line corresponds to $r=1.0$, the upper one to $r=4.0$.}
\label{Diagram}
\end{figure}
\end{center}

Considering pictures like this, one is faced with a striking 
qualitative similarity of the continuous
surfaces and their discrete counterparts.
Moreover, from numerical experiments it became clear that the
approximation is also quantitative: the
points on the discrete surface converge to the
corresponding points on the continous surface
as the mesh size $\epsilon$ goes to zero. 
The typical picture of the approximation
error is shown on Fig. \ref{Diagram}.

\epsfxsize=300pt

For a ``sufficiently typical'' family of discrete surfaces, 
we have plotted their maximal point distances from the respective 
continous counterpart. 
(See section \ref{Sect SG} for details.)
The slope of both curves is very close to one, which indicates
linear convergence.

All this suggests that
it might be possible to develop the classical differential geometry,
including both the theory of surfaces and of their transformations,
as a limit of the discrete constructions, just by refining the mesh size
in some directions. This is a common belief now, having however the status of
folklore only, since there are no rigorous mathematical statements supporting
it. On the other hand, the good quantitative properties of approximations
delivered by the discrete differential geometry suggest that they might
be put at the basis of the practical numerical algorithms for computations
in the differential geometry. Again, absence of mathematical results
on the quality of approximation prevents one from doing this.

The present paper aims at closing this gap. We take a step in this direction,
developing a numerical scheme for a class of nonlinear hyperbolic equations,
and proving general results on its convergence. It should be said that we
develop the numerics, having in mind also the above mentioned theoretical
applications to an alternative foundation of the differential geometry.
That means, first, that the class of hyperbolic systems we consider here
includes those coming from geometric applications, and, second, that our
discretizations respect the geometric structures, i.e. belong to the field
of discrete differential geometry. In particular, the notion of integrability,
in the guise of the multi--dimensional compatibility of discrete equations
\cite{BS}, will play a considerable role in our approach.

{\bf Example.} We shall illustrate our constructions by the well--known
Sine-Gordon equation:
\be\label{SineGordon}
\partial_x\partial_y\phi=\sin\phi.
\ee
A naive discretization of the Sine-Gordon equation could be obtained from
(\ref{SineGordon}) by replacing partial derivatives by their difference
analogs:
\be\label{dSGnaive}
\delta_x^{\eps}\delta_y^{\eps}\phi=\sin\phi,
\ee
where we introduced the following notation:
\be\label{deltas}
\delta^\eps_x p(x,y)=\frac{1}{\eps}\Big(p(x+\eps,y)-p(x,y)\Big),\quad
\delta^\eps_y p(x,y)=\frac{1}{\eps}\Big(p(x,y+\eps)-p(x,y)\Big).
\ee
In length, Eq. (\ref{dSGnaive}) reads:
\[
\phi(x+\eps,y+\eps)-\phi(x+\eps,y)-\phi(x,y+\eps)+\phi(x,y)=
        \eps^2\sin\phi(x,y).
\]
We shall prove an approximation theorem which implies that, on finite
domains, the solutions of a Goursat problem for (\ref{dSGnaive}) converge
with $\eps\to 0$ to the solutions of a Goursat problem for (\ref{SineGordon}),
provided the initial data on the characteristic lines converge. However,
this is not the whole story. Solutions of the Sine-Gordon equation correspond
to surfaces with constant negative Gaussian curvature. The discretization
(\ref{dSGnaive}) is non--geometric, and it not clear how to construct discrete
surfaces from its solutions. This is closely related to the fact that this
discretization does not inherit the integrability of the
Sine-Gordon equation. There exists a different one,
due to Hirota \cite{H}, which is itself a discrete integrable system. Its
geometric meaning and relation to discrete surfaces with the constant
negative Gaussian curvature was clarified by \cite{BP1}. The Hirota's
discretization of the Sine-Gordon equation reads:
\beqn
\lefteqn{\sin\frac{1}{4}\big(\phi(x+\eps,y+\eps)-\phi(x+\eps,y)-\phi(x,y+\eps)+
\phi(x,y)\big)} \nonumber\\
&=& \frac{\eps^2}{4}\sin\frac{1}{4}\big(\phi(x+\eps,y+\eps)+\phi(x+\eps,y)+
\phi(x,y+\eps)+\phi(x,y)\big).\qquad
\label{dSGHirota}
\eeqn
Our theory yields the convergence results for the solutions of this equation
to the solutions of (\ref{SineGordon}), which can be extended to the convergence
of surfaces, their associated families, and their B\"acklund transformations.
The existence of B\"acklund transformations is considered as a characteristic
property of integrability \cite{RS}.

It should be noticed that there exists an extensive literature dealing with
the numerical solution of the sine--Gordon and similar equations, see, e.g.
\cite{AHS1,AHS2,FV,FS,SV}. Most of these references have a physical background
and motivation, and due to this the problems settled and solved there are
different from those relevant to the differential geometry. In particular,
neither reference contains convergence results for the Goursat problem.
Notice, further, that our results are applicable to a large class
of hyperbolic systems describing various further geometries.

The structure of the paper is the following. In Sect. \ref{Sect 2d} we
formulate the continuous and discrete setup of the two--dimensional hyperbolic
systems and the corresponding Goursat problems. The $C^1$--convergence
result is proven which holds for all difference schemes with a local
approximation property. The $C^r$--approximation under the appropriate
conditions is established in Sect. \ref{Sect smooth}. The theory is extended to
the case of three independent variables in Sect. \ref{Sect 3D}. At this point
the notion of three--dimensional compatibility starts to play the key
role; it turns out to be intimately related to the integrability.
Therefore, the convergence result holds only for difference schemes with
these properties. The theory is illustrated in Sect. \ref{Sect SG}, where we
apply the convergence results to an integrable discretization of the
sine--Gordon equation, and thus prove the convergence of discrete $K$--surfaces
and their B\"acklund transformations to the continuous counterparts. Finally,
in Sect. \ref{Sect nD} the theory is extended to the case of an arbitrary
number of independent variables. The Appendix (Sect. \ref{Sect App})
contains the technical proofs of some statements of the main text.

Further pictures of discrete K-surfaces (like the ones
that constitute the data for Fig.\ref{Diagram}),
as well as a movie visualization of the convergence 
can be found on

{\tt http://www-sfb288.math.tu-berlin.de/$\sim$bobenko}$\;\;$.

\section{Two-dimensional theory}
\label{Sect 2d}

In this section we prove an approximation theorem for a certain class
of hyperbolic differential and difference equations in two dimensions.
Later on we will consider also more general $d$--dimensional systems.
The notations we will use for domains of independent variables
are the following: let ${\mathbf r}=(r_1,\ldots,r_d)$ consist of positive
numbers $r_i>0$, then
\be\label{Cdomain}
\Omega({\mathbf r})=[0,r_1]\times\ldots\times[0,r_d]\subset\rz^d.
\ee
As domains for discrete equations, we use parts of rectangular lattices
inside $\Omega({\mathbf r})$, with possibly different grid sizes along
different coordinate axes ${\bfEps}=(\eps_1,\ldots,\eps_d)$:
\be\label{Ddomain}
\Omega^{\bfeps}(\mathbf r)=[0,r_1]^{\eps_1}\times\ldots
\times[0,r_d]^{\eps_d}\subset\prod_{i=1}^d(\eps_i\gz).
\ee
where $[0,r]^{\eps}=[0,r]\cap(\eps\gz)$. The dependent variables of the
differential and difference equations under consideration are supposed to
belong to a normed vector space $\cX$, the norm in which will be denoted by
$|\cdot|$.

In the two--dimensional situation the notations will be, somewhat
inconsistently, simplified: we denote continuous domains by
$\Omega(r)=[0,r]\times[0,r]\subset\rz^2$, and discrete ones by
$\Omega^\eps(r)=[0,r]^{\eps}\times[0,r]^{\eps}\subset(\eps\gz)^2$,
Each $\Omega^\eps(r)$ contains $O({\eps}^{-2})$ of grid points. If
$\eps_2$ is an integer multiple of $\eps_1$, then
$\Omega^{\eps_2}(r)\subset\Omega^{\eps_1}(r)$. It will be convenient to assume
that $\eps$ attain only values of the form $2^{-k}$ with a positive integer $k$;
for such $\eps$, the relation $\eps_1<\eps_2$ implies that $\eps_2$ is an
integer multiple of $\eps_1$. We also define the limiting domains
\be
\Omega^0(r)=\bigcup_{\eps=2^{-k}} \Omega^\eps(r),
\ee
which are everywhere dense in $\Omega(r)$. Each point $x\in\Omega^0(r)$
belongs to $\Omega^\eps(r)$ with $\eps=2^{-k}$ for all $k$ large enough,
$k\ge k_0(x)$. This allows us to speak about convergence by $\eps=2^{-k}\to 0$
of sequences of functions $a^\eps$ defined on $\Omega^\eps(r)$. In case of
convergence, the limiting function may be thought of as defined on
$\Omega^0(r)$. If such a limiting function is continuous (Lipschitz)
on $\Omega^0(r)$, it can be extended to a continuous (resp. Lipschitz)
function on $\Omega(r)$.

Introduce the difference quotient operators $\delta^\eps_x$ and
$\delta^\eps_y$, acting on functions $p:\Omega^\eps(r)\rightarrow\cX$,
by the formulas (\ref{deltas}).
The functions $\delta^\eps_xp$ and $\delta^\eps_yp$ are
defined everywhere on $\Omega^\eps(r)$ except for the points with maximal
$x$ or $y$ coordinate.

\begin{dfn}
By a {\itbf continuous 2D hyperbolic system} we call a system of partial
differential equations for functions $a,b:\Omega(r)\rightarrow\cX$ of the form
\be\label{C2DH}
\partial_xa=f(a,b), \qquad \partial_yb=g(a,b),
\ee
with smooth enough functions $f,g:\cX\times\cX\rightarrow\cX$.
A {\itbf Goursat problem} consists of prescribing the initial values
\be\label{C2DI}
a(x,0)=a_0(x),\qquad b(0,y)=b_0(y)
\ee
for $x\in[0,r]$ and $y\in[0,r]$, respectively. The functions
$a_0,b_0:[0,r]\mapsto\cX$ are also supposed to be smooth enough.
\end{dfn}
\begin{dfn}
A {\itbf discrete 2D hyperbolic system} (or, better, a one--parameter family
of such systems) consists of two partial difference
equations for $a,b:\Omega^\eps(r)\rightarrow\cX$ of the form
\be\label{D2DH}
\delta^\eps_xa=f^{\eps}(a,b), \qquad \delta^\eps_yb=g^{\eps}(a,b),
\ee
with smooth functions $f^{\eps},g^{\eps}:\cX\times\cX\rightarrow\cX$.
A {\itbf Goursat problem} for this system consists of prescribing the initial
values
\be\label{D2DI}
a(x,0)=a_0^{\eps}(x),\qquad b(0,y)=b_0^{\eps}(y)
\ee
for $x\in[0,r]^{\eps}$ and $y\in[0,r]^{\eps}$, respectively.
\end{dfn}

{\bf Example.} Any equation of the type $\partial_x\partial_yu=F(u)$ can be
brought into the form (\ref{C2DH}) by a variety of substitutions
$a=A(u,\partial_xu)$ and $b=B(u,\partial_yu)$. For instance, a canonical way
to do this for the Sine-Gordon equation (\ref{SineGordon}) is to
introduce two new dependent variables
\be\label{AB}
a=\partial_x\phi, \qquad b=\phi,
\ee
which have to satisfy a system of the form (\ref{C2DH}):
\be\label{contABeqn}
\partial_ya=\sin b, \qquad \partial_xb=a.
\ee
The second order difference equation (\ref{dSGnaive}) can be dealt with in a
way which mimics the continuous situation: introduce two new dependent
variables
\be\label{dABnaive}
a=\delta_x^{\eps}\phi, \qquad b=\phi,
\ee
then they have to satisfy a discrete 2D hyperbolic system
\be\label{dABeqnnaive}
\delta_y^{\eps}a=\sin b, \qquad \delta_x^{\eps}b=a.
\ee
A similar procedure can be performed with the discretization
(\ref{dSGHirota}). Set
\beqn\label{dAB}
a=\delta_x^{\eps}\phi, \qquad b=\phi+\frac{\eps}{2}\,\delta_y^{\eps}\phi.
\eeqn
(The second equation may be written as
$b(x,y)=(\phi(x,y+\eps)+\phi(x,y))/2$.)
Then the following holds:
\beqn
b(x+\eps,y)-b(x,y) & = & \frac{\eps}{2}(a(x,y+\eps)+a(x,y)),\\
e^{\displaystyle i\eps a(x,y+\eps)/2}-
e^{\displaystyle i\eps a(x,y)/2} & = &
\frac{\eps^2}{4}\Big(e^{\displaystyle ib(x+\eps,y)}-
e^{\displaystyle -ib(x,y)}\Big).\qquad
\eeqn
(the first equation guarantees the existence of $\phi$ for a given pair of
functions $a$, $b$, satisfying (\ref{dAB}), while the second one is equivalent
to (\ref{dSGHirota})). These equations can be solved for $a(x,y+\eps)$ and
$b(x+\eps,y)$. The result reads:
\be\label{dABeqnHirota}
\delta_y^{\eps}a=\frac{2}{i\eps^2}\log\frac{1-(\eps^2/4)\exp(-ib-i\eps a/2)}
{1-(\eps^2/4)\exp(ib+i\eps a/2)},\qquad
\delta_x^{\eps}b=a+\frac{\eps}{2}\,\delta_y^{\eps}a.
\ee
Both discrete 2D hyperbolic systems (\ref{dABeqnnaive}), (\ref{dABeqnHirota})
approximate the continuous one (\ref{contABeqn}) in the sense of the next
definition.
\begin{dfn}
A discrete 2D hyperbolic system (\ref{D2DH}) {\itbf approximates} the
continuous one (\ref{C2DH}), if the functions $f^{\eps}$, $g^{\eps}$
continuously depend on $\eps\in(0,\eps_0)$, and
\be
f(a,b)=\lim_{\eps\to 0}f^{\eps}(a,b),\qquad
g(a,b)=\lim_{\eps\to 0}g^{\eps}(a,b),
\ee
uniformly on any compact subset of $\cX\times\cX$. If this convergence holds
in the $C^1$--topology, we speak about the $\mathbf C^1$--{\itbf approximation}. If
\be
f^\eps(a,b)=f(a,b)+\cO(\eps), \qquad g^\eps(a,b)=g(a,b)+\cO(\eps)
\ee
uniformly on any compact subset of $\cX\times\cX$, we speak about
the $\cO(\eps)$--{\itbf approximation}.
\end{dfn}

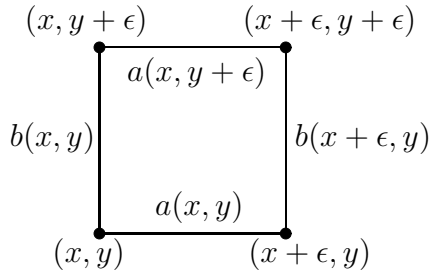
\begin{figure}[htbp]
\begin{center}
\setlength{\unitlength}{0.06em}
\begin{picture}(200,150)(-50,-20)
  \put(100,  0){\circle*{6}} \put(0  ,100){\circle*{6}}
  \put(  0,  0){\circle*{6}}  \put(100,100){\circle*{6}}
  \put( 0,  0){\line(1,0){100}}
  \put( 0,100){\line(1,0){100}}
  \put(  0, 0){\line(0,1){100}}
  \put(100, 0){\line(0,1){100}}
  \put(-25,-15){$(x,y)$}
  \put(80,-15){$(x+\eps,y)$}
  \put(77,110){$(x+\eps,y+\eps)$}
  \put(-40,110){$(x,y+\eps)$}
  \put(30,10){$a(x,y)$}
  \put(15,82){$a(x,y+\eps)$}
  \put(-48,47){$b(x,y)$}
  \put(105,47){$b(x+\eps,y)$}
\end{picture}
\caption{An elementary quadrilateral}\label{elem square}
\end{center}
\end{figure}

The following result is almost obvious:
\begin{prp}\label{exist2}
The Goursat problem for a discrete 2D hyperbolic system (\ref{D2DH})
has a unique solution $(a^\eps,b^\eps)$ on $\Omega^\eps(r)$.
\end{prp}
{\bf Proof.} At this point a useful remark should be done.
Although our notations might suggest that the variables $(a,b)$
are attached to the points of the two-dimensional lattice
$\Omega^\eps(r)$, it is more natural to assume that they are
attached to the {\em edges} of this lattice: $a(x,y)$ -- to the
horizontal edge connecting the vertices $(x,y)$ and $(x+\eps,y)$,
and $b(x,y)$ -- to the vertical edge connecting the vertices
$(x,y)$ and $(x,y+\eps)$. So, the equations (\ref{D2DH}) give the
fields sitting on the right and on the top edges of an elementary
square, provided the fields sitting on the left and on the bottom
ones are known. See Fig. \ref{elem square}. By induction, the
whole solution can be calculated, starting with the fields sitting
on all edges on the coordinate axes. \BBox
\vfree

Now the main result of this section can be formulated.
\begin{thm}
\label{mainthm}
Let a family of discrete 2D hyperbolic systems (\ref{D2DH})
$\cO(\eps)$--approximate the continuous 2D hyperbolic system (\ref{C2DH})
in the $C^1$ sense. Let also the discrete initial data (\ref{D2DI})
approximate the continuous ones (\ref{C2DI}):
\be
a_0^\eps(x)=a_0(x)+\cO(\eps),\qquad b_0^\eps(y)=b_0(y)+\cO(\eps)
\ee
uniformly for $x\in [0,r]^{\eps}$ and $y\in [0,r]^{\eps}$, respectively.
Then the sequence of solutions $(a^\eps,b^\eps)$ converges uniformly in
the following sense: there exist $\bar{r}\in(0,r]$ and a pair of
Lipschitz--continuous functions $(a,b)$ on $\Omega(\bar{r})$ such that
\be
a^\eps(x,y)=a(x,y)+\cO(\eps), \qquad
b^\eps(x,y)=b(x,y)+\cO(\eps)
\ee
for all $(x,y)\in\Omega^\eps(\bar{r})$. The functions $(a,b)$ solve the
continuous Goursat problem for (\ref{C2DH}) on $\Omega(\bar{r})$.
If the functions $f$, $g$ admit global Lipschitz constants, then one can
choose $\bar{r}=r$.
\end{thm}
{\bf Proof.}
In general one cannot expect $\bar{r}=r$ in the Theorem
because the solutions of the limiting equations may have blow-ups that
are absent in the discretization. Consequently, one essential step in
proving the theorem is to attain $\eps$-independent {\em \'a priori}
bounds on $a^\eps$ and $b^\eps$.

\begin{lem}\label{Nbound}
Let the norms of initial data $a^\eps_0$, $b^\eps_0$ be bounded by
$\eps$-independent constants. Then there exists $\bar{r}\in(0,r]$
such that the norms of the solutions $(a^\eps,b^\eps)$ are bounded
on the respective $\Omega^\eps(\bar{r})$ independently of $\eps$.
Furthermore, if $f$ and $g$ have a global Lipschitz constant $\cal L$
on the whole of $\cX\times\cX$, then one can choose $\bar{r}=r$.
\end{lem}
{\bf Proof of Lemma \ref{Nbound}.}
Let $|a_0^{\eps}|,|b_0^{\eps}|\le M_0$ with $M_0>0$.
We show that, fixing an arbitrary $P_0>M_0$, it is possible to find
$\bar{r}\leq r$ (independent of $\eps$) such that one has $|a^\eps(x,y)|,\,
|b^\eps(x,y)|\le P_0$ for $(x,y)\in\Omega^\eps(\bar{r})$. Actually, it is
enough to take
\be\label{rdef}
\bar{r}=\frac{P_0-M_0}{\cF(P_0)},
\ee
where
\be\label{cF}
\cF(M)=\sup_\eps\sup_{|a|,|b|<M}\{|f^\eps(a,b)|,|g^\eps(a,b)|\}.
\ee
Indeed, from the difference equations (\ref{D2DH}), written in length as
\beqn
a^\eps(x,y) & = & a^\eps(x,y-\eps)+
\eps f^\eps(a^\eps(x,y-\eps),b^\eps(x,y-\eps)), \label{aeqn}\\
b^\eps(x,y) & = & b^\eps(x-\eps,y)+
\eps g^\eps(a^\eps(x-\eps,y),b^\eps(x-\eps,y)),  \label{beqn}
\eeqn
we can conclude by induction that
\be\label{abest}
|a^\eps(x,y)|\le M_0+y\cF(P_0),\quad |b^\eps(x,y)|\le M_0+x\cF(P_0),
\ee
at least as long as the left--hand sides remain $\le P_0$. Since $M_0+\bar{r}
\cF(P_0)=P_0$, they remain $\le P_0$ for all $x,y\le\bar{r}$.

The argument for the case when $f$, $g$, and therefore $f^\eps$, $g^\eps$
have a global Lipschitz constant $\cL$, is a little bit different.
Set
\be
\Delta(x,y)=\max\{|a^\eps(x,y)|,|b^\eps(x,y)|,M_0\}.
\ee
From
(\ref{aeqn}), (\ref{beqn}) it follows readily that
\be
\Delta(x,y)\le \Big(1+\eps J(\Delta(x',y'))\Big)\cdot\Delta(x',y'),
\ee
where $(x',y')=(x,y-\eps)$ or $(x',y')=(x-\eps,y)$, depending on which of
$|a^\eps(x,y)|$ or $|b^\eps(x,y)|$ is greater, and where
\be\label{J(M)}
J(M)=\frac{1}{M}\sup_\eps\sup_{|a|,|b|\le M}\{|f^\eps(a,b)|,|g^\eps(a,b)|\}.
\ee
Now, this function admits for $M\ge M_0$ an estimate by an absolute constant:
\[
J(M)\leq\frac{1}{M_0}(|f^\eps(0,0)|+|g^\eps(0,0)|)+2\cL:=\cK.
\]
Hence,
\[
\Delta(x,y)\le (1+\eps\cK)\cdot\Delta(x',y').
\]
Using Lemma \ref{Gronwall} below (or a simple induction), we find:
\be
\Delta(x,y)\le M_0\exp(2\cK(x+y))\le M_0\exp(4\cK r).
\ee
This finishes the proof. \BBox
\vfree

On the last step we used the simple particular case $d=2$, $\kappa=0$ of
the following lemma, which will be used repeatedly later on. Its proof is
given in the Appendix.
\begin{lem}\label{Gronwall}
Let the function $\Delta:\Omega^{\bfeps}(\mathbf r)\mapsto\rz_+$ satisfy
the following condition: for any $x\in\Omega^{\bfeps}(\mathbf r)$ different
from the origin, there exists an index $1\le i\le d$ such that
\be
\Delta(x)\leq(1+\eps_i{\cal
K})\Delta(x-\eps_i\bfe_i)+\eps_i\kappa, \ee where $\cK,\kappa$ are
some non-negative numbers, and $\bfe_i$ is the unit vector of the
$i$-th coordinate axis. Then
\be
\Delta(x)\leq\max\Big(\Delta(0),\frac{\kappa}{\cal K}\Big)
\exp\Big(2{\cK}\sum_{j=1}^d x_j\Big).
\ee
\end{lem}

We demonstrate further that the solutions $(a^\eps,b^\eps)$ are not only
uniformly bounded, but are actually Lipschitz continuous with a
Lipschitz constant independent of $\eps$. In virtue of the equations
(\ref{D2DH}) and Lemma \ref{Nbound} it
is clear that the difference quotients $\delta^\eps_ya^\eps$ and
$\delta^\eps_xb^\eps$ are uniformly bounded. This turns out to be true
also for $\delta^\eps_xa^\eps$ and $\delta^\eps_yb^\eps$.
\begin{lem}\label{Lbound}
Let the initial data of the continuous Goursat problem
$a_0,b_0:[0,r]\mapsto\cX$ be $C^1$ functions, with the $C^1$--norm less then
$M>0$, and let the initial data
of the discrete Goursat problem $a_0^\eps,b_0^\eps:[0,r]^\eps\mapsto\cX$
satisfy
\be
|a_0^\eps(x)-a_0(x)|\leq M\eps,\quad |b_0^\eps(y)-b_0(y)|\leq M\eps.
\ee
With $\bar{r}\in(0,r]$ chosen according to Lemma \ref{Nbound},
the expressions $\delta^\eps_xa^\eps$ and $\delta^\eps_yb^\eps$
are bounded on the $\Omega^\eps(\bar{r})$ by $\eps$-independent constants.
\end{lem}
{\bf Proof.}
Let $M_1$ be a common bound on the values of the solutions of the
discrete Goursat problems.
Set
\be\label{M2}
M_2=\sup_\eps\sup_{|a|,|b|\le M_1}\Big\{|f^\eps(a,b)|,|g^\eps(a,b)|,
|\partial_af^\eps(a,b)|,\ldots, |\partial_bg^\eps(a,b)|\Big\}
\ee
(recall that $f^\eps\rightarrow f$ and $g^\eps\rightarrow g$ locally
uniformly in $C^1$). One can assume that $M$ from the condition of the
Lemma is greater than $M_1$ and $M_2$. By the mean value theorem, we find
an estimate for $\delta^\eps_xa^\eps(x,y)$ with $y=0$:
\be\label{3M}
|\delta_x^\eps a_0^\eps(x)| \leq |\delta_x^\eps a_0(x)|+
\eps^{-1} |a_0^\eps(x+\eps)-a_0(x+\eps)|
        +\eps^{-1}|a_0(x)-a_0^\eps(x)| \le 3M.
\ee
Proceeding from $y$ to $y+\eps$, we find:
\beq
|\delta^\eps_xa^\eps(x,y+\eps)|&\leq&|\delta^\eps_xa^\eps(x,y)|
                +\eps|\delta^\eps_xf^\eps(a^\eps(x,y),b^\eps(x,y))|\\
        &\leq&|\delta^\eps_xa^\eps(x,y)|+\eps M(|\delta^\eps_xa^\eps(x,y)|
                +|\delta^\eps_xb^\eps(x,y)|)\\
        &\leq&(1+\eps M)|\delta^\eps_xa^\eps(x,y)|
                +\eps M^2.
\eeq
Now Lemma \ref{Gronwall} yields the desired estimate:
\[
|\delta^\eps_xa^\eps(x,y)|\leq 4M\exp(M\bar{r}).
\]
The same reasoning applies to $\delta^\eps_yb^\eps$.  \BBox
\vfree

{\bf Proof of Theorem \ref{mainthm}, continued.}
We have to show that the $\eps$-dependent solutions of difference equations
have limits $(a^0,b^0)$ on $\Omega^0(\bar{r})$ with the convergence rates
\beqn
\label{uconv}
\sup_{(x,y)\in\Omega^\eps(\bar{r})}|a^\eps(x,y)-a^0(x,y)|,
|b^\eps(x,y)-b^0(x,y)|= \cO(\eps),
\eeqn
and that these limiting functions $(a^0,b^0)$ can be extended to continuous
functions $(a,b)$ on $\Omega(\bar{r})$ satisfying the differential equations.

Take $M>M_2$ (cf. (\ref{M2})) such that it bounds also $a^\eps$ and
$b^\eps$ along with their respective difference quotients (see Lemma
\ref{Lbound}). We will prove that $a^\eps(x,y)$ and $b^\eps(x,y)$ are
Cauchy sequences at any point $(x,y)\in\Omega^0(\bar{r})$.
To this end, fix $\eps,\eps'$ which are of the usual form $2^{-k}$;
we assume $\eps>\eps'$, so that $\eps/\eps'$ is a natural number. Below
we use the following abbreviations: $\Delta^{\eps\eps'}a(x,y)$
for $|a^\eps(x,y)-a^{\eps'}(x,y)|$, and $f^\eps[a^\eps,b^\eps](x,y)$ for
$f^\eps(a^\eps(x,y),b^\eps(x,y))$, etc. We have:
\beqn
\lefteqn{\Delta^{\eps\eps'}a(x,y) \leq } \nonumber\\
 &&\Delta^{\eps\eps'}a(x,y-\eps)+\Big|\eps f^\eps[a^\eps,b^\eps](x,y-\eps)-
     \eps'\sum_{\kappa=1}^{{\eps}/{\eps'}}
     f^{\eps'}[a^{\eps'},b^{\eps'}](x,y-\kappa\eps')\Big|=  \nonumber\\
 &&\Delta^{\eps\eps'}a(x,y-\eps)+  {\eps'}\sum_{\kappa=1}^{{\eps}/{\eps'}}
     \Big|f^\eps[a^\eps,b^\eps](x,y-\eps)-
     f^{\eps'}[a^{\eps'},b^{\eps'}](x,y-\kappa\eps')\Big|. \label{10}
\eeqn
For the terms in the last sum with a fixed $0<\kappa\leq{\eps}/{\eps'}$ we have:
\beqn
\lefteqn{|f^\eps[a^\eps,b^\eps](x,y-\eps)-
        f^{\eps'}[a^{\eps'},b^{\eps'}](x,y-\kappa\eps')|\leq}  \nonumber\\
 &   & |f^\eps[a^\eps,b^\eps](x,y-\eps)-
                f^{\eps'}[a^\eps,b^\eps](x,y-\eps)|  \label{l1}\\
 & + & |f^{\eps'}[a^\eps,b^\eps](x,y-\eps)-
                f^{\eps'}[a^{\eps'},b^{\eps'}](x,y-\eps)|  \label{l2}\\
 & + & |f^{\eps'}[a^{\eps'},b^{\eps'}](x,y-\eps)-
                f^{\eps'}[a^{\eps'},b^{\eps'}](x,y-\kappa\eps')|.  \label{l3}
\eeqn
Here $(\ref{l1})=\cO(\eps)$ because $f^\eps=f+\cO(\eps)$,
$(\ref{l3})=\cO(\eps)$ due to Lemma \ref{Lbound}, and
\[
(\ref{l2})\leq M\left(\Delta^{\eps\eps'}a(x,y-\eps)+
\Delta^{\eps\eps'}b(x,y-\eps)\right).
\]
Putting all this on the right--hand side of (\ref{10}), and taking into account
that the sum there contains $\cO(\eps/\eps')$ terms, we find:
\[
\Delta^{\eps\eps'}a(x,y)\leq (1+\eps M)\Delta^{\eps\eps'}a(x,y-\eps)+
\eps M\Delta^{\eps\eps'}b(x,y-\eps)+\cO(\eps^2).
\]
The analogous estimate holds for $b$:
\[
\Delta^{\eps\eps'}b(x,y)\leq (1+\eps M)\Delta^{\eps\eps'}b(x-\eps,y)+
\eps M\Delta^{\eps\eps'}a(x-\eps,y)+\cO(\eps^2).
\]
Introducing the quantities
\[
m_0:=\sup_{[0,\bar{r}]^\eps}\Big\{|a_0^\eps(\cdot)-a_0^{\eps'}(\cdot)|,
|b_0^\eps(\cdot)-b_0^{\eps'}(\cdot)|\Big\},
\]
\be\label{th1D1}
\Delta(x,y):=\max\{\Delta^{\eps\eps'}a(x,y),\Delta^{\eps\eps'}b(x,y),m_0\},
\ee
we get the estimate which holds also for the points of the boundary
of $\Omega^\eps(\bar{r})$ different from the origin:
\be\label{th1D2}
\Delta(x,y) \leq (1+\eps M)\Delta(x',y')+\cO(\eps^2),
\ee
where either $(x',y')=(x,y-\eps)$, or $(x',y')=(x-\eps,y)$. Now
Lemma \ref{Gronwall} yields the final estimate: $\Delta(x,y)=\cO(\eps)$.
Thus $a^\eps(x,y)$ and $b^\eps(x,y)$ are proved to form
Cauchy sequences at every point $(x,y)\in\Omega^0(\bar{r})$. Their limits
$a^0(x,y)$ and $b^0(x,y)$ satisfy (\ref{uconv}). These limiting functions
obviously have the same Lipschitz constants as all the $a^\eps$ and $b^\eps$.
Therefore there is a unique (Lipschitz) continuous extension of these
functions from $\Omega^0(\bar{r})$ to $\Omega(\bar{r})$; we call these
extensions $a$ and $b$.

It remains to be shown that $(a,b)$ solve the differential equations
(\ref{C2DH}). To do this, we prove that certain integral equations hold.
Let $(x,y)\in\Omega^0(\bar{r})$, then:
\beq
\lefteqn{a^\eps(x,y)=a_0^\eps(x)+\eps\sum_{k=0}^{y/\eps-1}
        f^\eps[a^\eps,b^\eps](x,k\eps)=}\\
  &  & a_0^\eps(x)+\eps\sum_{k=0}^{y/\eps-1}f[a,b](x,k\eps)+
  \eps\sum_{k=0}^{y/\eps-1}
        \Big(f^\eps[a^\eps,b^\eps](x,k\eps)-f[a,b](x,k\eps)\Big).
\eeq
All the terms inside the last sum are uniformly $\cO(\eps)$, therefore
the whole sum (with the pre-factor $\eps$) is also $\cO(\eps)$.
The first sum on the right--hand side is, up to $\cO(\eps)$, equal to
$\int_0^yf(a(x,\eta),b(x,\eta))d\eta$, which exists since $f$ is smooth and
$a,b$ are Lipschitz. Therefore,
\be\label{inteq a}
a(x,y)=a_0(x)+\int_0^y f[a,b](x,\eta)d\eta
\ee
on $\Omega^0(\bar{r})$. By continuity, this holds on all of $\Omega(\bar{r})$.
It follows that $a$ is everywhere differentiable with respect to $x$, and
$\partial_xa=f(a,b)$. The function $b$ is treated in the same manner.
\BBox

\begin{cor}
The two dimensional hyperbolic Goursat problem (\ref{C2DH}),
(\ref{C2DI}) possesses a unique classical solution.
\end{cor}
{\bf Proof.} The existence part of this statement is an immediate consequence
of Theorem \ref{mainthm}; uniqueness is easy to show as follows. Let $(a,b)$
and $(\hat{a},\hat{b})$ be two solutions of the system of integral equations
consisting of (\ref{inteq a}) and
\be\label{inteq b}
b(x,y)=b_0(y)+\int_0^x g[a,b](\xi,y)d\xi.
\ee
Subtracting the corresponding equations, we find after some simple
manipulations:
\[
\Delta(x,y)\leq L\left(\int_0^x \Delta(\xi,y)d\xi+
        \int_0^y \Delta(x,\eta)d\eta\right),
\]
where we introduced the deviation function
\beq
\Delta(x,y)=|a(x,y)-\hat{a}(x,y)|+|b(x,y)-\hat{b}(x,y)|,
\eeq
and $L$ is a common Lipschitz constant of $f,g$ over the range of values of
$a,b$. Now the 2D version of the classical Gronwall inequality
($d=2$ case of Lemma \ref{Gron} below) implies $\Delta(x,y)\equiv 0$. \BBox

\begin{lem}\label{Gron}
Let the continuous function $\Delta:\Omega({\mathbf r})\mapsto\rz$ satisfy
\be\label{GronwallX}
\Delta(x_1,\ldots,x_d)\leq L\sum_{j=1}^d\int_0^{x_j}
\Delta(x_1,\ldots,x_{j-1},\xi_j,x_{j+1},\ldots,x_d)d\xi_j+Q,
\ee
with some constants $L,Q\ge 0$. Then
\be
\Delta(x_1,\ldots,x_d)\leq 2Q\exp\Big(2dL\sum_{i=1}^d x_i\Big).
\ee
\end{lem}
Proof of this lemma is put in the Appendix.

\section{Additional Smoothness}
\label{Sect smooth}

In this section we show that the discrete techniques can be used
to prove regularity of the solutions to the hyperbolic equations.

 \begin{thm}\label{smoothm}
 Let the assumptions of Theorem \ref{mainthm} hold;
additionally, let the partial difference quotients up to the
order $k+1$ of the discrete initial data be uniformly bounded
independently of $\eps$:
\be
|(\delta^\eps_x)^ma_0^\eps(x)|,\;|(\delta^\eps_y)^nb_0^\eps(y)|\leq
M, \quad m,n\leq k+1, \ee
Suppose that the convergence $f^\eps\to f$, $g^\eps\to g$
is locally uniform in $C^{k+1}$. Then the limit functions
$a=\lim_{\eps\to 0}a^{\eps}$, $b=\lim_{\eps\to 0}b^{\eps}$ belong
to $C^k$, and moreover these limits are uniform in $C^k$:
\be
\sup_{\Omega^\eps(r-k\eps)}|(\delta^\eps_x)^m(\delta^\eps_y)^n
a^\eps-
        \partial_x^m\partial_y^m a|,\;
        |(\delta^\eps_x)^m(\delta^\eps_y)^n b^\eps-
        \partial_x^m\partial_y^m b|\to 0,\quad
        m+n\leq k.
\ee
\end{thm}

{\bf Remark.} Assume that the convergence $f^\eps\to f$, $g^\eps\to g$
is locally uniform in $C^\infty$ and that the initial data $a_0,b_0$ are
$C^{\infty}$--smooth. Then the canonical choice of the discrete initial
data $a_0^\eps(x)=a_0(x)$ and $b_0^\eps(y)=b_0(y)$ guarantees the
$C^k$--convergence for any $k$. One may then loosely speak of
$C^\infty$-approximation. \vfree

The ideas of the proof and even the essential estimates are
basically the same as one would most likely use in the continuous
setting. However, the discrete setting has the advantage that in
contrast to higher-order partial derivatives all the difference
quotients automatically exist, commute with each other etc.

First, one obtains {\it \'a priori} estimates not only for the
values of $a^\eps$, $b^\eps$, but for their higher order
difference quotients. We will need discrete analogues of the $C^k$-norms.
Let $\cY$ be a normed linear space; for a
function $u:\Omega^\eps(r)\rightarrow\cY$, define
 \be
 \nn u\nn_0= \sup_{\Omega^\eps(r)}|u|,
 \ee
 and
 \be\label{norms}
 \nn u\nn_K = \max_{k+\ell\leq K} \sup_{\Omega^\eps(r-K\eps)}
 |(\delta^\eps_x)^k(\delta^\eps_y)^\ell u|.
 \ee
 The following statement comes to replace the chain rule, which is
 no more available in the discrete context.
 \begin{lem}\label{composition}
 Let $f:\cX\times\cX\rightarrow\cX$ be a smooth function, and consider two functions
 $a,b:\Omega^\eps(r)\rightarrow\cX$. Then the $K$-th order
 difference quotients of $f(a,b)$ can be estimated as follows ($m+n=K$):
 \beqn \label{fest}
 |(\delta^\eps_x)^m(\delta^\eps_y)^nf[a,b](x,y)|&\leq&
        A\cdot\left(|(\delta^\eps_x)^m(\delta^\eps_y)^na(x,y)|+
        |(\delta^\eps_x)^m(\delta^\eps_y)^nb(x,y)|\right)+\nonumber \\
        &&P(\nn a\nn_{K-1},\nn b\nn_{K-1})+B\eps.
 \eeqn
 Here $P$ is a polynomial on two variables of total degree $\le K$
 with positive coefficients;
 the constants $A$, $B$, and the coefficients of $P$ depend only on $K$,
 on the Lipschitz constant of the functions $a$ and $b$, and on the expressions
 $\nn (D^k_aD^\ell_bf)[a,b]\nn_0$ for $k+\ell\le K+1$.
\end{lem}
The technical proof of the lemma in put in the Appendix. \vfree

{\bf Remark.} In applications of Lemma \ref{composition}, we have to handle
a whole $\eps$-dependent family of functions $a^\eps,b^\eps,f^\eps,g^\eps$
at once. Thus the constants $A$, $B$ and the polynomial $P$ become, in
principle, $\eps$-dependent, too. However, under the conditions like those
of Theorem \ref{smoothm}, one can choose these constants independently of
$\eps$. In particular, if
$\sup_\eps\nn a^\eps\nn_K$ and $\sup_\eps\nn b^\eps\nn_K$ are finite, then
$(\delta^\eps_x)^m(\delta^\eps_y)^nf^\eps(a^\eps,b^\eps)$ is bounded
for $m+n=K$, independently of $\eps$.

\begin{lem}
 Under the conditions of Theorem \ref{smoothm},
\be
\sup_\eps\nn a^\eps\nn_{k+1}<\infty,\hspace{10pt}
\sup_\eps\nn b^\eps\nn_{k+1}<\infty.
\ee
\end{lem}
 {\bf Proof} goes by induction with respect to the
 total degree $m+n=K\le k+1$ of the difference quotient.
 So assume that $\sup_\eps\nn a^\eps\nn_{K-1},
 \sup_\eps\nn b^\eps\nn_{K-1}<\infty$ is already proved.

 First, look at
 the difference quotients of $a^{\eps}$ with $n>0$. Then
 \be
 (\delta^\eps_x)^m(\delta^\eps_y)^na^\eps(x,y)=
        (\delta^\eps_x)^m(\delta^\eps_y)^{n-1}f^\eps(a^\eps,b^\eps)(x,y),
 \ee
 and the statement follows from the induction hypothesis by using Lemma
 \ref{composition}, as pointed out in the preceeding remark.
 Analogously, the difference quotients of $b^\eps$
 with $m>0$ are readily estimated.
 If $n=0$, then we have:
 \be
 (\delta^\eps_x)^ma^\eps(x,y)=(\delta^\eps_x)^ma^\eps(x,y-\eps)+
        \eps(\delta^\eps_x)^mf^\eps[a^\eps,b^\eps](x,y-\eps).
 \ee
 Apply Lemma \ref{composition} to find:
 \beq
 |(\delta^\eps_x)^ma^\eps(x,y)|&\leq&
    (1+\eps A)|(\delta^\eps_x)^ma^\eps(x,y-\eps)|\\
    &&+\eps A |(\delta^\eps_x)^mb^\eps(x,y)|\\
 &&+\eps P(\nn a^\eps\nn_{K-1},\nn b^\eps\nn_{K-1})+\eps^2B.
 \eeq
 By the induction hypotheses, the case considered above and the remark
 following Lemma \ref{composition}, the quantities on
 the right-hand side are controlled, so that Lemma
 \ref{Gronwall} can be applied, which yields a uniform,
 $\eps$-independent estimate for $(\delta^\eps_x)^ma^\eps$.
 Completely the same reasoning applies to $(\delta^\eps_y)^nb^\eps$.
 \BBox\vfree

 Theorem \ref{smoothm} is now a direct consequence of the next
 lemma.

\begin{lem}\label{HLbounds}
 Let a sequence of functions
 $\{u^\eps:\Omega^\eps(r)\rightarrow\cX\}_{\eps\in E_0}$
 be bounded in the discrete $C^{k+1}$-norm, independently of $\eps$:
 \be
 \sup_\eps\nn u^\eps\nn_{k+1}=M<\infty.
 \ee
 Then there
 exists a function $u\in C^k(\Omega(r),\cX)$, and a subsequence
 $\{u^\eps\}_{\eps\in E_\infty}$ for which
 \be
 \sup_{\Omega^{\eps}(r)}|(\delta^\eps_x)^m(\delta^\eps_y)^nu^\eps(x,y)-
        \partial_x^m\partial_y^nu(x,y)|\rightarrow 0\;\;
 \mbox{ as }\;\;\eps\rightarrow 0,\,\eps\in E_\infty.
 \ee
\end{lem}
{\bf Proof.}
The set $\Omega^0(r)$ consists of countably many points; choose any enumeration
\be
        \Omega^0(r)=\{(x_n,y_n)|n=1,2,3,\ldots\}.
\ee
Since all difference quotients up to the order $k$ are bounded at $(x_1,y_1)$,
there exists a subsequence $\{u^\eps\}_{\eps\in E_1}$ such that
\be
(\delta^\eps_x)^m(\delta^\eps_y)^n u^\eps(x_1,y_1)\rightarrow
u^{(m,n)}(x_1,y_1)\;\;\mbox{ as }\;\;\eps\rightarrow 0,\,\eps\in E_1.
\ee
for all $m+n\leq k$.
Repeat this procedure at $(x_2,y_2)$, where a convergent subsubsequence
with $E_2\subset E_1$ is selected, and so on. We get a series of infinite sets
$E_{i+1}\subset E_i$. Eventually, for all $m+n\le k$ the limits
$u^{(m,n)}(x,y)$ are defined everywhere on $\Omega^0(r)$.
Since  $\nn u^\eps\nn_{k+1}\le M$ holds $\eps$-uniformly,
the functions $u^{(m,n)}$ thus defined have the Lipschitz property on $\Omega^0(r)$,
so they possess unique continuous extensions to $\Omega(r)$ with $M$ as a
Lipschitz constant.

We are ready to describe the set $E_\infty$
from assertion of the lemma. We construct it as an infinite
sequence of numbers $\eps_j$ converging to zero, in such a way that
\be
\label{latticeconv}
\sup_{\Omega^{\eps_j}}
        |(\delta^{\eps_j}_x)^m(\delta^{\eps_j}_y)^n u^\eps(x,y)-u^{(m,n)}(x,y)|
        \le 2^{-j}(2M+1).
\ee
Let $\eps_0=1$, and take $\eps_j<\eps_{j-1}$ of the usual
form (integer power of $1/2$), satisfying
\be
|(\delta^{\eps_j}_x)^m(\delta^{\eps_j}_y)^nu^{\eps_j}(x',y')-
u^{(m,n)}(x',y')|\leq 2^{-j}
\ee
for all $m+n\le k$ at all points $(x',y')\in\Omega^{2^{-j}}(r)$
(we have with necessity $\eps_j\le 2^{-j}$, therefore $\Omega^{2^{-j}}(r)
\subset\Omega^{\eps_j}(r)$).
Such $\eps_j$ exists, because for finitely many lattice sites of
$\Omega^{2^{-j}}$ there holds a pointwise convergence argument from the
beginning of this proof.
In order to establish the estimate (\ref{latticeconv}) on the larger set
$\Omega^{\eps_j}(r)$, we make use of the Lipschitz properties.
Let $(x,y)\in\Omega^{\eps_j}(r)$ be arbitrary; then there always exists a point
$(x',y')\in\Omega^{2^{-j}}(r)$ with $|x'-x|+|y'-y|\le 2^{-j}$. Thus we conclude
\beq
\lefteqn{|(\delta^{\eps_j}_x)^m(\delta^{\eps_j}_y)^n u^{\eps_j}(x,y)
-u^{(m,n)}(x,y)|}\\
      &\le  & |(\delta^{\eps_j}_x)^m(\delta^{\eps_j}_y)^n
u^{\eps_j}(x,y)-(\delta^{\eps_j}_x)^m(\delta^{\eps_j}_y)^n u^{\eps_j}(x',y')|\\
       & + & |(\delta^{\eps_j}_x)^m(\delta^{\eps_j}_y)^nu^{\eps_j}(x',y')-
       u^{(m,n)}(x',y')|\\
       &+ & |u^{(m,n)}(x',y')-u^{(m,n)}(x,y)|\\
        &\le& M2^{-j}+2^{-j}+M2^{-j}.
\eeq
It remains to show that $u^{(m,n)}$ are indeed the respective partial
derivatives of $u$. This can be done by proving the corresponding integral
equations: from
\be
(\delta^\eps_x)^m(\delta^\eps_y)^nu^\eps(x,y)=
        (\delta^\eps_x)^m(\delta^\eps_y)^nu^\eps(0,y)+
        \eps\sum_{k=0}^{x/\eps-1}(\delta^\eps_x)^{m+1}(\delta^\eps_y)^n
    u^\eps(k\eps,y)
\ee
there follows in the limit $\eps\rightarrow 0, \eps\in E_\infty$
\be
u^{(m,n)}(x,y)=
        u^{(m,n)}(0,y)+\int_0^xu^{(m+1,n)}(\xi,y)d\xi,
\ee
so that $\partial_xu^{(m,n)}(x,y)=u^{(m+1,n)}(x,y)$.
\BBox\vfree

{\bf Proof of Theorem \ref{smoothm}.} With Lemma \ref{HLbounds} at hand,
the convergence of the difference quotients of $a^\eps$, $b^\eps$ is proved as
follows. Any subsequence of $(a^{\eps},b^{\eps})$ has a subsubsequence
which converges uniformly with all its difference quotients to a $C^k$
limit function. But by Theorem \ref{mainthm}, this limit function is
the unique solution $(a,b)$ of the correspondent continuous Goursat problem.
If any subsequence of a given sequence has a subsubsequence converging to
one and the same limit, then the whole sequence must converge to that
same limit.
\BBox

\section{Three--dimensional theory: \\
approximating B\"acklund transformations}
\label{Sect 3D}

The Sine-Gordon equation (\ref{SineGordon}) is a very special one, in that
it posseses B\"acklund transformations. This means that from a given
solution $\phi$ we can construct new solutions by solving ordinary
differential equations only. The famous formula for a family of elementary
B\"acklund transformations $\phi\mapsto\widetilde{\phi}$ for
(\ref{SineGordon}) reads:
\be\label{SGBT}
\partial_x\widetilde{\phi}+\partial_x\phi=
2\alpha\sin\frac{\widetilde{\phi}-\phi}{2},\qquad
\partial_y\widetilde{\phi}-\partial_y\phi=
\frac{2}{\alpha}\sin\frac{\widetilde{\phi}+\phi}{2}.
\ee
A direct calculation shows that this system is compatible (i.e.
$\partial_y(\partial_x\widetilde{\phi})=
\partial_x(\partial_y\widetilde{\phi})$), provided $\phi$ is a solution of
the Sine-Gordon equation, and then $\widetilde{\phi}$ is also a solution. An
equivalent way to express this state of affairs is to introduce, along with
the variables $a,b$ from (\ref{AB}), also the auxilary function
$\theta=(\widetilde{\phi}-\phi)/2$, which satisfies the
following system of ordinary differential equations:
\be\label{BT}
\partial_x\theta=-a+\alpha\sin\theta, \qquad
\partial_y\theta=\frac{1}{\alpha}\sin(b+\theta).
\ee
This system is compatible in the sense that
$\partial_y(\partial_x\theta)=\partial_x(\partial_y\theta)$, provided $(a,b)$
solve the system (\ref{contABeqn}) equivalent to the Sine--Gordon equation,
and the initial value
\be
\theta(0,0)=\theta_0
\ee
defined uniquely a solution of (\ref{BT}). Then the formulas
\be\label{BTfin}
\widetilde{a}=a+2\partial_x\theta=-a+2\alpha\sin\theta,\qquad
\widetilde{b}=b+2\theta
\ee
deliver a new solution $(\widetilde{a},\widetilde{b})$ of the 2D hyperbolic
system (\ref{contABeqn}) equivalent to the Sine--Gordon equation. Clearly,
B\"acklund transformations can be iterated in a straightforward manner.

\begin{dfn}
A {\itbf continuous 2D hyperbolic system with a B\"acklund tranformation}
is a compatible system of partial differential and difference equations
\beqn
\partial_ya=f(a,b), & \quad & \partial_xb=g(a,b),  \nonumber\\
\partial_x\theta=u(a,\theta), & \quad & \partial_y\theta=v(b,\theta),
\label{C3DH}\\
\delta_za=\xi(a,\theta), & \quad & \delta_zb=\eta(b,\theta)\nonumber
\eeqn
for functions $a,b,\theta:\Omega(r,R)\mapsto\cX$, where
\be
\Omega(r,R)=\{(x,y,z)\;|\;(x,y)\in\Omega(r),\;z=0,1,\ldots,R\}.
\ee
Here $f,g,u,v,\xi,\eta:\cX\times\cX\mapsto\cX$ are asuumed to be smooth
functions. A {\itbf Goursat problem} is posed by the requirement
\be\label{C3DI}
a(x,0,0)=a_0(x),\;\; b(0,y,0)=b_0(y),\;\; \theta(0,0,z)=\theta_0(z)
\ee
for $x\in[0,r]$, $y\in[0,r]$, and $z\in\{0,1,\ldots,R\}$, respectively,
with given smooth functions $a_0(x)$, $b_0(y)$ and a sequence
$\theta_0(0),\ldots,\theta_0(R-1)$.
\end{dfn}
The {\itbf compatibility conditions} mentioned in this definition,
are to assure the existence of solutions of the above Goursat problem.
They follow from
\[
\partial_y(\partial_x\theta)=\partial_x(\partial_y\theta),\quad
\partial_y(\delta_za)=\delta_z(\partial_ya),\quad
\partial_x(\delta_zb)=\delta_z(\partial_xb).
\]
In length, these conditions for (\ref{C3DH}) to be B\"acklund transformations
of the 2D hyperbolic system read:
\beqn
\lefteqn{D_a u(a,\theta)\cdot f(a,b)+D_\theta u(a,\theta)\cdot v(b,\theta)=}
\nonumber\\
 & = & D_b v(b,\theta)\cdot g(a,b)+D_\theta v(b,\theta)\cdot u(a,\theta),
\nonumber\\
\lefteqn{D_a\xi(a,\theta)\cdot f(a,b)+D_\theta\xi(a,\theta)\cdot v(b,\theta)=}
\nonumber\\
 & = & f\Big(a+\xi(a,\theta),b+\eta(b,\theta)\Big)-f(a,b), \label{C3DC}\\
\lefteqn{D_b\eta(b,\theta)\cdot g(a,b)+D_\theta\eta(b,\theta)\cdot u(a,\theta)=}
\nonumber\\
 & = & g\Big(a+\xi(a,\theta),b+\eta(b,\theta)\Big)-g(a,b).\nonumber
\eeqn
The existence of B\"acklund transformations may be regarded as one of the
possible definitions of the {\it integrability} of a given 2D continuous
hyperbolic system. For a given 2D continuous hyperbolic system with B\"acklund
transformations, not every discretization possesses the analogous property.
For instance, the naive discretization (\ref{dSGnaive}) of the Sine-Gordon
equation does not admit B\"acklund transformations, while the integrable
discretization (\ref{dSGHirota}) does. The difference analogs of the formulas
(\ref{SGBT}) read:
\beqn
\lefteqn{\sin\frac{1}{4}
\big(\widetilde{\phi}(x+\eps,y)-\widetilde{\phi}(x,y)+\phi(x+\eps,y)-
\phi(x,y)\big)=} \nonumber\\
&=& \frac{\eps\alpha}{2}\sin\frac{1}{4}
\big(\widetilde{\phi}(x+\eps,y)+\widetilde{\phi}(x,y)-
\phi(x+\eps,y)-\phi(x,y)\big).\qquad
\label{BTdSGx}\\
\lefteqn{\sin\frac{1}{4}
\big(\widetilde{\phi}(x,y+\eps)-\widetilde{\phi}(x,y)-\phi(x,y+\eps)+
\phi(x,y)\big)=} \nonumber\\
&=& \frac{\eps}{2\alpha}\sin\frac{1}{4}
\big(\widetilde{\phi}(x,y+\eps)+\widetilde{\phi}(x,y)+
\phi(x,y+\eps)+\phi(x,y)\big).\qquad
\label{BTdSGy}
\eeqn
Obviously, in the limit $\eps\to 0$ these equations approximate (\ref{SGBT}).
A very remarkable feature is that these equations closely resemble the
original difference equation (\ref{dSGHirota}), if one considers the tilde
as the shift in the third $z$--direction. Upon introducing the quantity
$\theta=(\widetilde{\phi}-\phi)/2$, one rewrites (\ref{BTdSGx}), (\ref{BTdSGy})
in the form of the system of first order equations approximating
(\ref{BT}), (\ref{BTfin}):
\beqn
\delta_x^{\eps}\theta & = & -a+\frac{1}{i\eps}
\log\frac{1-(\eps\alpha/2)\exp(-i\theta+i\eps a/2)}
{1-(\eps\alpha/2)\exp(i\theta-i\eps a/2)},   \label{dBTx}\\
\delta_y^{\eps}\theta & = & \frac{1}{i\eps}
\log\frac{1-(\eps/2\alpha)\exp(-ib-i\theta)}
{1-(\eps/2\alpha)\exp(ib+i\theta)},   \label{dBTy}
\eeqn
and
\be\label{dBTfin}
\widetilde{a}=a+2\delta_x^\eps\theta,  \qquad
\widetilde{b}=b+2\theta+\eps\delta_y^\eps\theta.
\ee
This suggests the following definition.
\begin{dfn}
A {\itbf discrete 3D hyperbolic system} is a collection of compatible
partial difference equations of the form
\beqn
\nonumber
\delta^\eps_ya=f^\eps(a,b), & \quad & \delta^\eps_xb=g^\eps(a,b),\\
\label{D3DH}
\delta^\eps_x\theta=u^\eps(a,\theta), & \quad &
\delta^\eps_y\theta=v^\eps(b,\theta),\\
\nonumber
\delta_za=\xi^\eps(a,\theta), & \quad & \delta_zb=\eta^\eps(b,\theta),
\eeqn
for functions $a,b,\theta:\Omega^\eps(r,R)\mapsto\cX$, where
\be
\Omega^\eps(r,R)=\{(x,y,z)\;|\;(x,y)\in\Omega^\eps(r),\;z=0,1,\ldots,R\}.
\ee
Here the functions $f^\eps,g^\eps,u^\eps,v^\eps,\xi^\eps,\eta^\eps:
\cX\times\cX\rightarrow\cX$ are smooth enough. A {\itbf Goursat problem}
consists of prescribing the initial data
\beqn
\label{D3DI}
a(x,0,0)=a_0^\eps(x),\;\;b(0,y,0)=b_0^\eps(y),\;\;\theta(0,0,z)=\theta_0^\eps(z)
\eeqn
for $x\in[0,r]^\eps$, $y\in[0,r]^\eps$, and $z\in\{0,1,\ldots,R\}$,
respectively.
\end{dfn}

{\itbf Compatibility conditions} are necessary for solutions of (\ref{D3DH})
to exist. These conditions express the following identities that have to be
fulfilled for the solutions:
\[
\delta_y^\eps(\delta_x^\eps\theta)=\delta^\eps_x(\delta^\eps_y\theta),\quad
\delta_y^\eps(\delta_za)=\delta_z(\delta_y^\eps a),\quad
\delta_x^\eps(\delta_zb)=\delta_z(\delta_x^\eps b).
\]
In length, these formulas read:
\beqn
\lefteqn{u^\eps\Big(a+\eps f^\eps(a,b),\theta+\eps v^\eps(b,\theta)\Big)-
u^\eps(a,\theta)=}
\nonumber\\
& = & v^\eps\Big(b+\eps g^\eps(a,b),\theta+\eps u^\eps(a,\theta)\Big)-
v^\eps(b,\theta),\nonumber\\
\lefteqn{\xi^\eps\Big(a+\eps f^\eps(a,b),\theta+\eps v^\eps(b,\theta)\Big)-
\xi^\eps(a,\theta)=}\nonumber\\
& = & \eps f^\eps\Big(a+\xi^\eps(a,\theta),b+\eta^\eps(b,\theta)\Big)-
\eps f^\eps(a,b),
\label{D3DC}\\
\lefteqn{\eta^\eps\Big(b+\eps g^\eps(a,b),\theta+\eps u^\eps(a,\theta)\Big)-
\eta^\eps(b,\theta)=}\nonumber\\
& = & \eps g^\eps\Big(a+\xi^\eps(a,\theta),b+\eta^\eps(b,\theta)\Big)-
\eps g^\eps(a,b).
\nonumber
\eeqn
This has to be satisfied identically in $a,b,\theta\in\cX$.

As demonstrated in \cite{BS}, the compatibility of a discrete 3D hyperbolic
system is closely related to its integrability in the sense of the soliton
theory. Moreover, such a key attribute of integrability as a {\it discrete
zero curvature representation with a spectral parameter} can be derived
from the fact of compatibility.

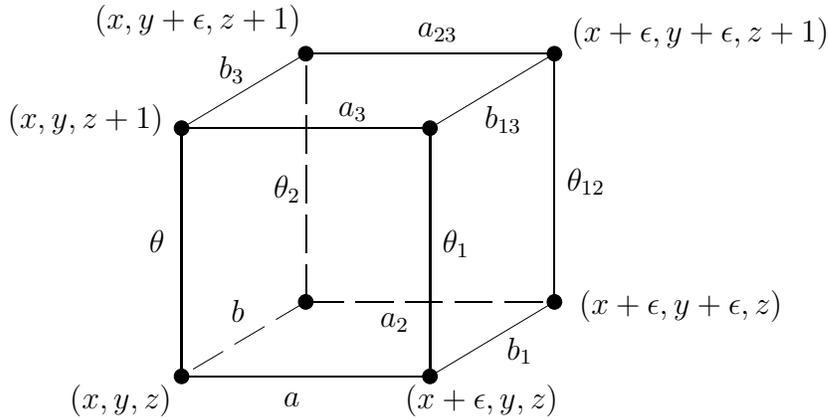
\begin{figure}[htbp]
\begin{center}
\setlength{\unitlength}{0.08em}
\begin{picture}(200,190)(-30,-20)
  \put(100,  0){\circle*{6}} \put(0  ,100){\circle*{6}}
  \put( 50, 30){\circle*{6}} \put(150,130){\circle*{6}}
  \put(  0,  0){\circle*{6}}  \put(100,100){\circle*{6}}
  \put( 50,130){\circle*{6}}  \put(150, 30){\circle*{6}}
  \put( 0,  0){\line(1,0){100}}
  \put( 0,100){\line(1,0){100}}
  \put(50,130){\line(1,0){100}}
  \multiput(50,30)(20,0){5}{\line(1,0){15}}
  \put(  0, 0){\line(0,1){100}}
  \put(100, 0){\line(0,1){100}}
  \put(150,30){\line(0,1){100}}
  \multiput(50,30)(0,20){5}{\line(0,1){15}}
  \put(  0,100){\line(5,3){50}}
  \put(100,100){\line(5,3){50}}
  \put(100,  0){\line(5,3){50}}
  \multiput(50,30)(-16.67,-10){3}{\line(-5,-3){12}}
     \put(-45,-13){$(x,y,z)$}
     \put(90,-13){$(x+\eps,y,z)$}
     \put(-70,100){$(x,y,z+1)$}
     \put(160,25){$(x+\eps,y+\eps,z)$}
     \put(-35,140){$(x,y+\eps,z+1)$}
     \put(157,135){$(x+\eps,y+\eps,z+1)$}
     \put(41,-12){$a$}
     \put(80,20){$a_2$}
     \put(63,105){$a_3$}
     \put(95,137){$a_{23}$}
     \put(-13,50){$\theta$}
     \put(37,72){$\theta_2$}
     \put(105,50){$\theta_1$}
     \put(155,75){$\theta_{12}$}
     \put(20,22){$b$}
     \put(131,7){$b_1$}
     \put(15,120){$b_3$}
     \put(122,100){$b_{13}$}
\end{picture}
\caption{Three-dimensional consistency}\label{fig.cube}
\end{center}
\end{figure}

\begin{prp}\label{exist3}
The Goursat problem for a discrete 3D hyperbolic system (\ref{D3DH})
satisfying the compatibility conditions (\ref{D3DC}) has a unique solution
$(a^\eps,b^\eps,\theta^\eps)$ on $\Omega^\eps(r,R)$.
\end{prp}
{\bf Proof.} Like in the proof of Proposition \ref{exist2}, it is enough
to demonstrate that the solution can be propagated along an elementary ``cube''
of the three--dimensional lattice, then the induction will end the proof. Again,
it is convenient to assume that the variables $a(x,y,z)$, $b(x,y,z)$,
$\theta(x,y,z)$ are attached not to the points $(x,y,z)\in\Omega^\eps(r,R)$,
but rather to the edges $[(x,y,z),(x+\eps,y,z)]$, $[(x,y,z),(x,y+\eps,z)]$,
$[(x,y,z),(x,y,z+1)]$, respectively. Denote (in this proof only) shifts
of the edge variables in the directions of $x,y,z$ axes by the
subscripts $1,2,3$, respectively. (See Fig. \ref{fig.cube}.) Then
the values $(a_2,b_1)$ are detrmined by the first equation in
(\ref{D3DH}), the values $(\theta_1,\theta_2)$ -- by the second one, and
the values $(a_3,b_3)$ -- by the third one. For the each one of the values
$(a_{23},b_{13},\theta_{12})$ sitting on the edges adjacent to
$(x+\eps,y+\eps,z+1)$ we get two possible values (from two equations defined
on two facets sharing the corresponding edge). The compatibility conditions
(\ref{D3DC}) guarantee that these two values for each of the edge variables
actually coincide. \BBox

\begin{thm}\label{3Dthm}
Let the family of discrete 3D hyperbolic systems (\ref{D3DH}) satisfying the
compatibility conditions (\ref{D3DC}) approximate the continuous 2D hyperbolic
system with a B\"acklund transformation (\ref{C3DH}). Let the approximation
be of the order $\cO(\eps)$ in the $C^1$--sense, so that the relations
\be
f^\eps(a,b)=f(a,b)+\cO(\eps),\;\;\mbox{ etc.}
\ee
and similar relations for the first partial derivatives hold uniformly on
any bounded set of $(a,b,\theta)\in\cX^3$. Let also
the initial data be approximated according to
\beq
a_0^\eps(x)=a_0(x)+\cO(\eps),\;\;\;b_0^\eps(y)=b_0(y)+\cO(\eps),\;\;\;
        \theta_0^\eps(z)=\theta_0(z)+\cO(\eps)
\eeq
uniformly for $x\in[0,r]^\eps$, $y\in[0,r]^\eps$, $z\in\{0,1,\ldots,R\}$,
respectively.

Then, for some $\bar{r}\in(0,r]$, the sequence of solutions
$(a^\eps,b^\eps,\theta^\eps)$ has a uniform limit
of Lipschitz--continuous functions $(a,b,\theta)$ on $\Omega(r,R)$
in the sense that the relations
\beq
a^\eps(x,y,z)&=&a(x,y,z)+\cO(\eps),\\
b^\eps(x,y,z)&=&b(x,y,z)+\cO(\eps),\\
\theta^\eps(x,y,z)&=&\theta(x,y,z)+\cO(\eps)
\eeq
hold uniformly in $(x,y,z)$ from the respective $\Omega^\eps(\rho,R)$.
Furthermore, $(a,b,\theta)$ solve the Goursat problem for the
continuous 2D system with a sequence of B\"acklund transformations.
\end{thm}
{\bf Proof} parallels the proof of Theorem \ref{mainthm}, and starts with
an {\em \'a priori} estimate for $a^\eps$, $b^\eps$, $\theta^\eps$.
\begin{lem}\label{Nbound2}
Let the norms of initial data $a^\eps_0$, $b^\eps_0$, $\theta^\eps_0$ be
bounded by $\eps$-independent constants. Then there exists $\bar{r}\in(0,r]$
such that the norms of the solutions $(a^\eps,b^\eps,\theta^\eps)$ are bounded
on the respective $\Omega^\eps(\bar{r},R)$ independently of $\eps$. If the
right--hand sides of Eqs. (\ref{D3DH}) possess a global Lipschitz constant
$\cal L$ on the whole of $\cX\times\cX\times\cX$, then one can choose
$\bar{r}=r$.
\end{lem}
{\bf Proof of Lemma \ref{Nbound2}.}
Let $|a_0^{\eps}|,|b_0^{\eps}|,|\theta^\eps_0|\le M_0$ with $M_0>0$. Set
\be
\cF(M)=\sup_\eps\sup_{|a|,|b|,|\theta|<M}\{|f^\eps(a,b)|,\ldots,
|\eta^\eps(b,\theta)|\}.
\ee
Choose an arbitrary $P_0>M_0$, and define inductively $P_{j+1}=P_j+\cF(P_j)$
for $j=0,1,\ldots,R-1$. Finally, set
\be\label{rdef2}
\bar{r}=\min_{j=0,1,\ldots,R}\frac{P_j-M_0}{2\cF(P_j)},
\ee
so that
\be\label{rprop2}
M_0+2\bar{r}\cF(P_j)\le P_j\quad {\rm for\;\;all}\quad j=0,1,\ldots,R.
\ee
We show now by induction in $z$ that for all $(x,y,z)\in\Omega^\eps(r,R)$
there holds:
\be\label{induct2}
|a^{\eps}(x,y,z)|,\,|b^{\eps}(x,y,z)|,\,|\theta^\eps(x,y,z)|\le P_z.
\ee
Indeed, like in Lemma \ref{Nbound} (see (\ref{abest})), from
the first pair of equations in (\ref{D3DH}) we conclude that
\be\label{abest2}
|a^\eps(x,y,0)|\le M_0+y\cF(P_0),\quad |b^\eps(x,y,0)|\le M_0+x\cF(P_0),
\ee
and similarly from the second pair of equations in (\ref{D3DH})
we conclude that
\be\label{thetest2}
|\theta^\eps(x,y,0)|\le M_0+(x+y)\cF(P_0),
\ee
as long as the left--hand sides of these inequalities remain $\le P_0$.
Due to (\ref{rprop2}), they remain $\le P_0$ for all
$(x,y)\in\Omega^\eps(\bar{r})$. Now assuming (\ref{induct2}) for a given $z$,
the third pair of equations in (\ref{D3DH}) immediately implies:
\be
|a^\eps(x,y,z+1)|,\,|b^\eps(x,y,z+1)|\le P_z+\cF(P_z)\le P_{z+1},
\ee
while from the second pair we derive:
\be
|\theta^\eps(x,y,z+1)|\le M_0+(x+y)\cF(P_{z+1})\le P_{z+1}.
\ee
This proves (\ref{induct2}), and therefore the first statement of lemma.
The second one is proved similarly to the analogous statement of Lemma
\ref{Nbound}. \BBox\vfree

Further, like in the two--dimensional case, we have to demonstrate that
not only the solutions $(a^\eps,b^\eps,\theta^\eps)$ are bounded on
$\Omega^\eps(\bar{r},R)$, but also their difference quotients do.
\begin{lem}\label{Lbound2}
Let the initial data of the continuous Goursat problem
$a_0,b_0:[0,r]\mapsto\cX$ be $C^1$ functions, with the $C^1$--norm less then
$M>0$, and let the initial data of the discrete Goursat problem
$a_0^\eps,b_0^\eps:[0,r]^\eps\mapsto\cX$ satisfy
\be
|a_0^\eps(x)-a_0(x)|\leq M\eps,\quad |b_0^\eps(y)-b_0(y)|\leq M\eps.
\ee
Let $\bar{r}\in(0,r]$ be chosen according to Lemma \ref{Nbound2}. Then
all the difference quotients $\delta^\eps_x a^\eps$, $\delta^\eps_y a^\eps$,
$\delta^\eps_x b^\eps$, $\delta^\eps_y b^\eps$, $\delta^\eps_x \theta^\eps$,
and $\delta^\eps_y \theta^\eps$ are bounded on all of the
$\Omega^\eps(\bar{r},R)$ by $\eps$-independent constants.
\end{lem}
{\bf Proof of Lemma \ref{Lbound2}.} Let the estimates (\ref{induct2}) hold
for all $(x,y,z)\in\Omega^\eps(\bar{r},R)$. Then from the equations
(\ref{D3DH}) it is immediately seen that
\[
|\delta^\eps_y a^\eps|,\,|\delta^\eps_x b^\eps|,\,
|\delta^\eps_x \theta^\eps|,\,|\delta^\eps_y \theta^\eps|\le \cF(P_z).
\]
Therefore, like in Lemma \ref{Lbound}, we only have to estimate the
quotients $\delta^\eps_x a^\eps$ and $\delta^\eps_y b^\eps$. Let $M$ be
chosen to be greater than the common Lipschitz constant of the functions
$f^\eps$, $g^\eps$, $\xi^\eps$ and $\eta^\eps$ for all $(a,b,\theta)$ with
$|a|,|b|,|\theta|\le P_R$. Then, proceeding as in the proof of Lemma
\ref{Lbound}, we have:
\beq
\lefteqn{|\delta^\eps_xa^\eps(x,y,0)|\leq|\delta^\eps_xa^\eps(x,y-\eps,0)|
                +\eps|\delta^\eps_xf^\eps[a^\eps,b^\eps](x,y-\eps,0)|}\\
        &\leq&|\delta^\eps_xa^\eps(x,y-\eps,0)|+\eps
    M(|\delta^\eps_xa^\eps(x,y-\eps,0)|
                +|\delta^\eps_xb^\eps(x,y-\eps,0)|)\\
        &\leq&(1+\eps M)|\delta^\eps_xa^\eps(x,y-\eps,0)|
                +\eps M\cF(P_0).
\eeq
By (\ref{3M}) and Lemma \ref{Gronwall}, we obtain:
\be
|\delta^\eps_x a^\eps(x,y,0)|\leq(3M+\cF(P_0))e^{2My}.
\ee
Proceeding from $z-1$ to $z$, we find:
\beq
\lefteqn{|\delta^\eps_xa^\eps(x,y,z)|\leq|\delta^\eps_xa^\eps(x,y,z-1)|
                +|\delta^\eps_x\xi^\eps[a^\eps,\theta^\eps](x,y,z-1)|}\\
        &\leq&|\delta^\eps_xa^\eps(x,y,z-1)|+
    M(|\delta^\eps_xa^\eps(x,y,z-1)|
                +|\delta^\eps_x\theta^\eps(x,y,z-1)|)\\
        &\leq&(1+M)|\delta^\eps_xa^\eps(x,y,z-1)|
                +M\cF(P_{z-1}).
\eeq
Applying again  Lemma \ref{Gronwall}, we find:
\be
|\delta^\eps_x a^\eps(x,y,z)|\leq
(|\delta^\eps_x a^\eps(x,y,0)|+\cF(P_{R-1}))e^{2Mz}
\leq (3M+2P_R)e^{2M(y+z)}.
\ee
An estimate for $|\delta^\eps_y b^\eps(x,y,z)|$ is obtained in a completely
similar way. \BBox\vfree

{\bf Proof of Theorem \ref{3Dthm}, continued.} Proceeding parallel to the
proof of the two--dimensional Theorem \ref{mainthm}, we introduce the
quantities $\Delta^{\eps\eps'}a=|a^\eps-a^{\eps'}|$ etc., and
\[
\Delta(x,y,z)=\max\{\Delta^{\eps\eps'}a(x,y,z),\Delta^{\eps\eps'}b(x,y,z),
\Delta^{\eps\eps'}\theta(x,y,z),m_0\},
\]
where $m_0$ is the supremum of the analogous quantities
$\Delta^{\eps\eps'}a_0=|a_0^\eps-a_0^{\eps'}|$ for the initial data (cf.
(\ref{th1D1})). The following inequality comes to replace (\ref{th1D2}):
\be\label{th3D2}
\Delta(x,y,z) \leq (1+\delta M)\Delta(x',y',z')+\cO(\eps\delta),
\ee
where either $(x',y',z')=(x-\eps,y.z)$ or $(x',y',z')=(x,y-\eps.z)$, in which
cases $\delta=\eps$, or else  $(x',y',z')=(x,y.z-1)$, and then $\delta=1$.
In any case Lemma \ref{Gronwall} is applicable and yields the final estimate:
$\Delta(x,y,z)=\cO(\eps)$. This proves that $a^\eps(x,y,z)$, $b^\eps(x,y,z)$
and $\theta^\eps(x,y,z)$ form Cauchy sequences at every point
$(x,y,z)\in\Omega^0(\bar{r},R)$. The end of the proof is completely analogous
to that of Theorem \ref{mainthm}. \BBox\vfree

\section{Approximation theorems for K-Surfaces}
\label{Sect SG}

In the present section we apply the theory developed so far to prove
that the known construction of discrete surfaces of constant negative
Gauss curvature $K=-1$ ({\itbf K--surfaces}, for short) may be actually used
not only to modelling the geometric properties of their continuous
counterparts, but also to a quantitative approximation. First we briefly
recall the correspondent geometric notions.

Let $F$ be a K-surface parametrized by its asymtotic lines:
\be
F:\Omega(r)\mapsto\rz^3.
\ee
This means that the vectors $\partial_x F$, $\partial_y F$, $\partial_x^2 F$,
$\partial_y^2 F$ are orthogonal to the normal vector $N:\Omega(r)\mapsto
S^2$. Reparametrizing the asymptotic lines, if necessary, we assume that
$|\partial_x F|=1$ and $|\partial_y F|=1$. The angle $\phi=\phi(x,y)$ between
the vectors $\partial_x F$, and $\partial_y F$ satisfies
the sine--Gordon equation (\ref{SineGordon}). Moreovere,
a K-surface is determined by a solution to (\ref{SineGordon}) essentially
uniquely. The correspondent construction is as follows. Consider the matrices
$U,V$ defined by the formulas
\beqn
U(a;\lambda) & = & \frac{i}{2}\left(\begin{array}{cc} a & -\lambda \\
-\lambda & -a \end{array}\right), \label{contLaxU}\\
V(b;\lambda) & = & \frac{i}{2}\left(\begin{array}{cc} 0 &
\lambda^{-1}\exp(ib)\\
        \lambda^{-1}\exp(-ib) & 0 \end{array}\right), \label{contLaxV}
\eeqn
taking values in the twisted loop algebra
\[
g[\lambda]=\{\xi:\rz_*\mapsto {\rm su}(2):\;
\xi(-\lambda)=\sigma_3\xi(\lambda)\sigma_3\},\quad
\sigma_3=\left(\begin{array}{cc} 1 & 0 \\ 0 & -1 \end{array}\right).
\]
Suppose now that $a$ and $b$ are some real--valued functions on $\Omega(r)$.
Then the zero curvature condition
\be\label{ZCC}
\partial_yU-\partial_xV+[U,V]=0
\ee
is satisfied identically in $\lambda$, if and only if $(a,b)$ satisfy
the system (\ref{contABeqn}), or, in other words, if $a=\partial_x\phi$ and
$b=\phi$, where $\phi$ is a solution of (\ref{SineGordon}). Given a solution
$\phi$, that is, a pair of matrices (\ref{contLaxU}), (\ref{contLaxV})
satisfying (\ref{ZCC}), the following system of linear differential equations
is uniquely solvable:
\be\label{Frame}
\partial_x\Phi=U\Phi,\quad   \partial_y\Phi=V\Phi,\quad
\Phi(0,0,\lambda)={\bf 1}.
\ee
Here $\Phi:\Omega(r)\mapsto G[\lambda]$ takes values in the twisted loop
group
\[
G[\lambda]=\{\Xi:\rz_*\mapsto {\rm SU}(2):\;
\Xi(-\lambda)=\sigma_3\Xi(\lambda)\sigma_3\}.
\]
The solution $\Phi(x,y;\lambda)$ yields the immersion $F(x,y)$ by the
{\itbf Sym formula}:
\be\label{contSym}
F(x,y)=\left(2\lambda\Phi(x,y;\lambda)^{-1}\partial_\lambda\Phi(x,y;\lambda)
\right)\Big|_{\lambda=1}.
\ee
(Here the canonical identification of ${\rm su}(2)$ with $\rz^3$ is used.)
Moreover, the right--hand side of (\ref{contSym}) at the values of
$\lambda$ different from $\lambda=1$ delivers a whole family of immersions
$F_\lambda:\Omega(r)\mapsto\rz^3$, all of which turn out to be
asymptotic lines parametrized K--surfaces. These surfaces $F_\lambda$
constitute the so--called {\itbf associated family} of $F$. \vfree

Now we turn to {\itbf discrete K-surfaces}. Let $F^\eps$ be a discrete
surface parametrized by asymptotic lines, i.e. an immersion
\be
F^\eps:\Omega^\eps(r)\rightarrow\rz^3
\ee
such that for each $(x,y)\in \Omega^\eps(r)$ the five points
$F^\eps(x,y)$, $F^\eps(x\pm\eps,y)$, and $F^\eps(x,y\pm\eps)$ lie in a single
plane $\cP(x,y)$. It is required that all edges of the discrete surface
$F^\eps$ have the same length $\eps\ell$, that is
$|\delta_x^\eps F^\eps|=|\delta_y^\eps F^\eps|=\ell$, and it turns out to be
convenient to assume that $\ell=(1+\eps^2/4)^{-1}$.
The same relation we presented between K-surfaces and solutions to the
(classical) sine--Gordon equation (\ref{SineGordon}) can be found
between discrete K-surfaces and solutions to the sine--Gordon equation
in Hirota's discretization (\ref{dSGHirota}): Consider the matrices
$\cU^\eps$, $\cV^\eps$ defined by the formulas
\beqn
\cU^\eps(a;\lambda) & = & (1+\eps^2\lambda^2/4)^{-1/2}\left(\begin{array}{cc}
        \exp(i\eps a/2) &-i\eps\lambda/2 \\
        -i\eps\lambda/2 & \exp(-i\eps a/2) \end{array}\right),
    \label{discrU}\\
\cV^\eps(b;\lambda) & = & \!(1+\eps^2\lambda^{-2}/4)^{-1/2}\left(\begin{array}{cc}
        1 & (i\eps\lambda^{-1}/2)\exp(ib) \\
    (i\eps\lambda^{-1}/2)\exp(-ib) & 1 \end{array}\right)\!.\nonumber\\
    \label{discrV}
\eeqn
Let $a$, $b$ be real--valued functions on $\Omega^\eps(r)$, and consider the
discrete zero curvature condition
\be\label{dZCC}
\cU^\eps(x,y+\eps;\lambda)\cdot\cV^\eps(x,y;\lambda)=
\cV^\eps(x+\eps,y;\lambda)\cdot\cU^\eps(x,y;\lambda)
\ee
(where $\cU^\eps$ and $\cV^\eps$ depend on $(x,y)\in\Omega^\eps(r)$
through the dependence of $a$ and $b$ on $(x,y)$, respectively).
A direct calculation shows that (\ref{dZCC}) is equivalent to the system
(\ref{dABeqnHirota}), or, in other words,
to the equation (\ref{dSGHirota}) for the function $\phi$ defined by
(\ref{dAB}). The formula (\ref{dZCC}) is the compatibility condition of
the following system of linear difference equations:
\beqn\label{discFrame}
\Psi^\eps(x+\eps,y;\lambda) & = & \cU^\eps(x,y;\lambda)\Psi^\eps(x,y;\lambda),
\nonumber\\
\Psi^\eps(x,y+\eps;\lambda) & = & \cV^\eps(x,y;\lambda)\Psi^\eps(x,y;\lambda),\\
\Psi^\eps(0,0;\lambda) & = & {\bf 1}.\nonumber
\eeqn
So, any solution of (\ref{dSGHirota}) uniquely defines a matrix
$\Psi^\eps:\Omega^\eps(r)\mapsto G[\lambda]$ satisfying (\ref{discFrame}).
This can be used to finally construct the immersion
by an analog of the Sym formula:
\be\label{discrSym}
F^\eps(x,y)=\left(2\lambda\Psi^\eps(x,y;\lambda)^{-1}
\partial_\lambda\Psi^\eps(x,y;\lambda)\right)\Big|_{\lambda=1}.
\ee
The geometric meaning of the function $\phi$ is the following: the angle
between the edges $F^\eps(x+\eps,y)-F^\eps(x,y)$ and
$F^\eps(x,y+\eps)-F^\eps(x,y)$
is equal to $(\phi(x+\eps,y)+\phi(x,y+\eps))/2$; the angle between the edges
$F^\eps(x,y+\eps)-F^\eps(x,y)$ and $F^\eps(x-\eps,y)-F^\eps(x,y)$ is equal to
$\pi-(\phi(x,y+\eps)+\phi(x-\eps,y))/2$; the angle between the edges
$F^\eps(x-\eps,y)-F^\eps(x,y)$ and $F^\eps(x,y-\eps)-F^\eps(x,y)$ is equal to
$(\phi(x-\eps,y)+\phi(x,y-\eps))/2$; and the angle between the edges
$F^\eps(x,y-\eps)-F^\eps(x,y)$ and $F^\eps(x+\eps,y)-F^\eps(x,y)$ is equal to
$\pi-(\phi(x,y-\eps)+\phi(x+\eps,y))/2$. In particular, these angles sum up to
$2\pi$, so that the four neighboring vertices of $F^\eps(x,y)$ lie in one
plane, as they should. Again, the right--hand side of (\ref{discrSym})
at the values of $\lambda$ different from $\lambda=1$ delivers an
associated family $F_\lambda^\eps$ of discrete
asymptotic lines parametrized K--surfaces. \vfree

Finally, we discuss the {\itbf B\"acklund transformations} for continuous
and discrete K--surfaces. Introduce the matrix
\be\label{DBack}
\cW(\theta;\lambda)=\left(\begin{array}{cc}
\alpha\exp(i\theta) & -i\lambda \\ -i\lambda & \alpha\exp(-i\theta)
\end{array}\right).
\ee
It is easy to see that the matrix differential equations
\be
\partial_x\cW=\widetilde{U}\cW-\cW U,\qquad
\partial_y\cW=\widetilde{V}\cW-\cW V
\ee
are equivalent to the formulas (\ref{BT}), (\ref{BTfin}). On the other hand,
these matrix differential equations constitute a sufficient condition for
the solvability of the system consisting of (\ref{Frame}) and
\be
\widetilde{\Phi}=\cW\Phi.
\ee
So, frames $\Phi$ can be in a consistent way extended into the third
direction $z$ (shift in which is encoded by the tilde), which results
also in the transformation of the K--surfaces $F\mapsto\widetilde{F}$,
and moreover of the whole associated family, via (\ref{contSym}).

Similarly, the matrix equations
\beqn
\cW(x+\eps,y;\lambda)\cU^\eps(x,y;\lambda) & = &
\widetilde{\cU}^\eps(x,y;\lambda)\cW(x,y;\lambda),\\
\cW(x,y+\eps;\lambda)\cV^\eps(x,y;\lambda) & = &
\widetilde{\cV}^\eps(x,y;\lambda)\cW(x,y;\lambda)
\eeqn
are equivalent to the formulas (\ref{dBTx}), (\ref{dBTy}), (\ref{dBTfin}),
and, on the other hand, assure the solvability of the system consisting of
(\ref{discFrame}) and
\be
\widetilde{\Psi}^\eps=\cW\Psi^\eps.
\ee
Therefore, also the frames $\Psi^\eps$ of the discrete surfaces can be
extended in the third direction $z$. This leads to the transformation
of discrete K--surfaces and their associated families, according to
(\ref{discrSym}).

Geometrical meaning of B\"acklund transformations is the following (see, e.g.,
\cite{BP1}). The asymptotic lines parametrizations of $F$ and $\widetilde{F}$
correspond, and the vector $\Delta F=\widetilde{F}-F$ lies in the intersection
of tangential planes of $F$ and of $\widetilde{F}$ in the corresponding points.
This vector $\Delta F$ has a constant length $\alpha$, and $\theta$ is the
angle between this vector and one of the asymptotic directions. Classically,
a B\"acklund transformation is completely defined by $\alpha$ and the value of
$\theta$ in one point. The definitions in the discrete case are completely
analogous.
\vfree

Now we are prepared to state the approximation theorem for K-surfaces.

\begin{thm}\label{Kapprox}
Let $a_0:[0,r]\mapsto \rz$ and $b_0:[0,r]\mapsto S^1=\rz/(2\pi\gz)$ be two
smooth functions. Then:
\begin{itemize}
\item There exists a unique asymptotic line parametrized K--surface
$F:\Omega(r)\rightarrow\rz^3$ such that its characteristic angle
$\phi: \Omega(r)\mapsto S^1$ on the coordinate axes satisfies:
\be
\partial_x\phi(x,0)=a_0(x),\quad \phi(0,y)=b_0(y),\quad x,y\in[0,r].
\ee
\item For any $\eps>0$ there exists a unique asymptotic line parametrized
discrete K--surface $F^\eps:\Omega^\eps(r)\rightarrow\rz^3$ such that its
characteristic angle $\phi^\eps: \Omega^\eps(r)\mapsto S^1$ on the coordinate
axes satisfies:
\be
\label{discreteinitialdata}
\phi^\eps(x+\eps,0)-\phi^\eps(x,0)=\eps a_0(x),\quad
\phi^\eps(0,y+\eps)+\phi^\eps(0,y)=2b_0(y),
\ee
for $ x,y\in[0,r-\eps]^\eps$.
\item There holds:
\be\label{KdKapprox}
\sup_{\Omega^\eps(r)}|F^\eps-F|\le C\eps,
\ee
where $C$ does not depend on $\eps$. Moreover, for pair $(m,n)$ of
nonnegative integers, there holds
\be\label{KdKapproxCk}
\sup_{\Omega^\eps(r-k\eps)}|(\delta^\eps_x)^m(\delta^\eps_y)^n F^\eps
-\partial_x^m\partial_y^n F|\rightarrow 0\quad {\rm as}\quad \eps\to 0.
\ee
\item The estimates (\ref{KdKapprox}), (\ref{KdKapproxCk}) still hold,
uniformly for $\lambda\in[\Lambda^{-1},\Lambda]$ with any $\Lambda>1$,
if one replaces in these estimates the immersions $F$, $F^\eps$ by their
associated families $F_\lambda$, $F^\eps_\lambda$, respectively.
\item Finally, if $\theta_1,\ldots,\theta_R\in S^1$ are parameters of a
sequence of B\"acklund transformations, then for the transformed surfaces
there also hold estimates analogous to (\ref{KdKapprox}), (\ref{KdKapproxCk}).
\end{itemize}
\end{thm}

{\bf Remark:}
Local uniform convergence of the discrete surfaces to the
continous still holds if we replace $a_0$ and $b_0$ in
(\ref{discreteinitialdata}) by some $\eps$-close initial
data. However, we lose the $C^\infty$-convergence of the $F^\eps$;
it is weakened to local uniform convergence of the frames in the
matrix norm.
Replacing the discrete initial data by sequences
$a^\eps_0$, $b^\eps_0$ that converge in the discrete $C^k$-sense
to $a_0$, $b_0$ leads to local $C^k$-convergence of $F^\eps$ to the
limit; this is a direct consequence of the following proof in combination
with Theorem \ref{smoothm}.
\vfree

Eventually, we can be explicit about what is shown in
Fig.\ref{Diagram}. The family of K-surfaces under consideration is in
correspondence to solutions of the discrete Sine-Gordon-Equation for the
same intial data $a_0$, $b_0$ on different domains $\Omega^\eps(r)$. 
The lower curve corresponds to $r=1$, and the upper curve to $r=4$.
The mesh size was chosen $\eps=2^{-k}$ for $k=5,6,\ldots,11$. 
As initial data we prescribed two smooth functions with no (apparent)
special properties.
The most precise approximation $F^{\eps_*}$ with $\eps_*\approx 0.002$
is used as the reference point (the limiting smooth surface).
The error in Fig. \ref{Diagram} is the $C^0$ distance of the
map $F^\eps$ to $F^{\eps_*}$.

For special classes of K-surfaces possessing additional symmetries
the convergence may be much faster. In particular, the rate of
converegence for the Amsler surface in Fig.\ref{Amsler} is apparently
quadratic in $\epsilon$. 
This ``superconvergence'' might be explained
by the fact that the discretization preserves the defining
geometric (surface contains two straight asymptotic lines) and
analytic (for a relation to the Painlev\'e III equation, see
\cite{PIII}) properties.
\vfree

{\bf Proof.} Theorems \ref{mainthm} and \ref{smoothm}
yield the existence and the uniqueness of solutions $(a^\eps,b^\eps)$ to
the difference equations on the whole of $\Omega^\eps(r)$, the existence
and uniqueness of the solutions $(a,b)$ to the differential equations,
and the $C^\infty$ approximation of the latter by the former.

It remains to prove that similar approximation holds also for the immersions
$F^\eps$, $F$. To do this, we prove the approximation property for the frames
$\Psi^\eps$, $\Phi$, uniformly in $\lambda\in[\Lambda^{-1},\Lambda]$, and
then use the Sym formula. Recall that these frames are defined as solutions
of the Cauchy problems for the system of linear difference equations
(\ref{discFrame}), respectively for the system of linear differential
equations (\ref{Frame}). Since the zero curvature conditions (\ref{dZCC}),
(\ref{ZCC}) are satisfied, the existence of $\Psi^\eps$, $\Phi$ is guaranteed
by standard ODE theory. Furthermore, at any point $(x,y)$,
$\Psi^\eps(\lambda)$ and $\Phi(\lambda)$ are analytic functions of
$\lambda\in D$, where $D$ is a closed disc in the complex plane of $\lambda$
that contains all points having distance $(2\Lambda)^{-1}$ from the interval
$[\Lambda^{-1},\Lambda]$. Since $D$ has positive distance to $i\rz$ and $\infty$,
all the matrices $\cU^\eps$, $\cV^\eps$, and $U$, $V$ are bounded uniformly
with respect to $\lambda\in D$ and $(x,y)\in\Omega(r)$. It is easy to see that
\be
\cU^\eps(a;\lambda)={\bf 1}+\eps U(a;\lambda)+\cO(\eps^2),\quad
\cV^\eps(b;\lambda)={\bf 1}+\eps V(b;\lambda)+\cO(\eps^2).
\ee
Even more is true: If $a^\eps=a+\cO(\eps)$ and $b^\eps=b+\cO(\eps)$, then
\be\label{Laxapprox}
\cU^\eps(a^\eps;\lambda)={\bf 1}+\eps U(a;\lambda)+\cO(\eps^2),\quad
\cV^\eps(b^\eps;\lambda)={\bf 1}+\eps V(b;\lambda)+\cO(\eps^2).
\ee
To estimate $\Psi^\eps-\Phi$, observe first that
\[
\Phi(x+\eps,y)=\Phi(x,y)+\int_x^{x+\eps}U(\xi,y)\Phi(\xi,y)d\xi=
       ({\bf 1}+\eps U(x,y))\Phi(x,y)+\cO(\eps^2).
\]
On the other hand, due to (\ref{Laxapprox}),
\[
\Psi^\eps(x+\eps,y)=\cU^\eps(x,y)\Psi^\eps(x)=
({\bf 1}+\eps U(x,y))\Psi^\eps(x,y)+\cO(\eps^2).
\]
Therefore,
\be
\Psi^\eps(x+\eps,y)-\Phi(x+\eps,y)=({\bf 1}+\eps U(x,y))(\Psi^\eps(x,y)-
\Phi(x,y))+\cO(\eps^2).
\ee
Similarly,
\be
\Psi^\eps(x,y+\eps)-\Phi(x,y+\eps)=({\bf 1}+\eps V(x,y))(\Psi^\eps(x,y)-
\Phi(x,y))+\cO(\eps^2).
\ee
Due to the zero curvature condition (\ref{ZCC}), one can find $\Phi(x,y)$
by first integrating the first equation in (\ref{Frame}) from $(0,0)$ to $(x,0)$
along the $x$--axis, and then integrating the second equation
in (\ref{Frame}) from $(x,0)$ to $(x,y)$ parallel to the $y$--axis.
The similar holds for $\Psi^\eps$. Therefore the classical Gronwall inequality
can be used to conclude that
\be
\Psi^\eps(x,y)-\Phi(x,y)=\cO(\eps).
\ee
For real $\lambda$, the frames $\Psi^\eps$ and $\Phi$ belong to
${\rm SU}(2)$, so it is immediate that also
\be
(\Psi^\eps)^{-1}(x,y)-\Phi^{-1}(x,y)=\cO(\eps).
\ee
Furthermore, since $\Phi(\lambda)$ und $\Psi^\eps(\lambda)$ are analytic
functions of $\lambda\in D$, we have for any $(x,y)\in \Omega^\eps(r)$ by
Cauchy's theorem:
\beqn
\lefteqn{\sup_{\lambda\in[\Lambda^{-1},\Lambda]}
\nn\partial_\lambda\Psi^\eps(x,y;\lambda)-\partial_\lambda\Phi(x,y;\lambda)\nn
\leq }\nonumber\\
 & & 2\Lambda\sup_{\lambda\in D}\nn\Psi^\eps(x,y;\lambda)-\Phi(x,y;\lambda)\nn
=\cO(\eps).
\eeqn
The last two estimates imply that, for all $\lambda\in[\Lambda^{-1},\Lambda]$
and uniformly on the respective $\Omega^\eps(r)$,
\be
F^\eps_\lambda-F_\lambda=2\lambda(\Psi^\eps(\lambda))^{-1}
\partial_\lambda\Psi^\eps(\lambda)-
        2\lambda(\Phi(\lambda))^{-1}\partial_\lambda\Phi(\lambda)=\cO(\eps).
\ee
It remains to prove the approximation of the higher order partial derivatives
of $F$ by the correspondent difference quotients of $F^\eps$. Introducing the
notations
\[
{\bf U}^\eps=(\cU^\eps-{\bf 1})/\eps=U+\cO(\eps),\quad
{\bf V}^\eps=(\cV^\eps-{\bf 1})/\eps=V+\cO(\eps),
\]
it is easy to see that for corresponding solutions $(a^\eps,b^\eps)$ and
$(a,b)$ one has discrete $C^k$-approximation for all $k>0$ and
all $\lambda\in D$:
\[
\nn {\bf U}^\eps-U\nn_k\rightarrow 0\mbox{ and }
\nn {\bf V}^\eps-V\nn_k\rightarrow 0
\]
as $\eps\rightarrow 0$.
We find for $m+n=k+1$ with $m>0$:
\beq\label{frameinduction}
|(\delta^\eps_x)^m(\delta^\eps_y)^n(\Psi^\eps-\Phi)|=
   |(\delta^\eps_x)^{m-1}(\delta^\eps_y)^n({\bf U}^\eps\Psi^\eps-U\Phi)|
   +\cO(\eps)\qquad\nonumber\\
\leq C\nn U\nn_k\cdot\nn\Psi^\eps-\Phi\nn_k+
    C\nn {\bf U}^\eps-U\nn_k\cdot\nn\Psi^\eps\nn_k+\cO(\eps).
\eeq
Here we used that for discrete $C^k$-norms of matrix products
\be
\nn A\cdot B\nn_k\le C_k\nn A\nn_k\cdot\nn B\nn_k
\ee
holds (cf. the remark after the proof of Lemma \ref{composition} in the Appendix).
If $m=0$, we can do the same calculations with the roles of $x$ and $y$ interchanged
and $V$, ${\bf V}^\eps$ in place of $U$, ${\bf U}^\eps$.
>From this estimate, we conclude by induction in $k$ that
\[
\nn\Psi^\eps-\Phi\nn_k\rightarrow 0,
\]
and therefore
\[
\lim_{\eps\to 0}\sup_{\Omega^\eps(r)}
|(\delta^\eps_x)^m(\delta^\eps_y)^n \Psi^\eps-
\partial_x^m\partial_y^n \Phi|=0.
\]
Again by the Cauchy estimate, we get also the similar result for the
respective $\lambda$-derivatives for all values
$\lambda\in[\Lambda^{-1},\Lambda]$. From the Sym formulas, we get
(\ref{KdKapproxCk}). Finally, the statement about the approximation
of B\"acklund transformed surfaces follows in a completely similar way
with the refence to Theorem \ref{3Dthm}. \BBox
\vfree

As a corollary, we automatically get the non-trivial classical theorem
on the permutability of B\"acklund transformations, due to Bianchi: given
a K--surface $F$ and its two B\"acklund transformations $F_1$, $F_2$, there
exists a unique K--surface $F_{12}$ which is a B\"acklund transformation of
$F_1$ and of $F_2$. This follows now by a continuous limit in two directions
of a four--dimensional compatible lattice system.

\section{General hyperbolic systems}
\label{Sect nD}

The theory developped so far can be generalized to higher dimensions without
difficulties. We want here formulate the most general setting in which our
techniques can be used to prove the approximation results for
discretizations of Goursat problems for nonlinear hyperbolic systems.
As an illustrative example, the reader can think of the equation
\[
\partial_x\partial_y\partial_z u=F(u,\partial_x u,\partial_y u,\partial_z u,
\partial_x\partial_y u,\partial_x\partial_z u,\partial_y\partial_z u).
\]
The Goursat problem for it consists of prescribing the values of
\[
u(x,y,0),\;u(x,0,z),\;u(0,y,z)\quad {\rm for}\quad 0\le x,y,z\le r.
\]
The above equation can be rewritten as a hyperbolic system:
\[
\left\{\begin{array}{l}
\partial_x u =a, \quad \partial_y u =b, \quad \partial_z u =c,\\
\partial_y a =h, \quad \partial_z b =f, \quad \partial_x c =g,\\
\partial_z a =g, \quad \partial_x b =h, \quad \partial_y c =f,\\
\partial_xf =\partial_y g =\partial_z h =F(u,a,b,c,f,g,h).
\end{array}\right.
\]
A natural (naive) way to discretize it consists of replacing all partial
derivatives $\partial_x$ etc. by the correspondent difference quotients
$\delta_x^\eps$ etc. In the so obtained hyperbolic difference system it is
natural to assume that the variables $a,b,c$ live on the edges of the
cubic lattice starting from the point $(x,y,z)$ in the direction of the
axes $x,y,z$, respectively, and that the variables $f,g,h$ are defined
on the two--cells (elementary squares) adjacent to the point $(x,y,z)$
and orthogonal to the axes $x,y,z$, respectively. The generalization of this
construction is as follows.

We use the notations (\ref{Cdomain}) for the domain of the differential
system, and (\ref{Ddomain}) for the domain of its discretization. We denote
by ${\bf e}_i$ the vector whose only nonvanishing component is 1 in the
$i$th position. The dependent variables are denoted by
$\vec{a}=(a_1,\ldots,a_N)\in\cX_1\times\ldots\cX_N$.
\begin{dfn}
For each $1\le k\le N$, let there be chosen a nonempty subset
$\cE_k\subset\{1,\ldots,d\}$ of independent variables; let $\cD_k$ be its
complement, so that there is a disjoint union
\be\label{division}
\{1,\ldots,d\}=\cE_k\cup \cD_k.
\ee
A {\itbf discrete d--dimensional hyperbolic system} is a collection of
compatible difference equations
\be\label{GDDH}
\delta^{\eps_i}_{x_i}a_k=f_{(k,i)}(\vec{a}), \quad i\in\cE_k,
\ee
for the functions $a_k:\Omega^{\bfeps}(\mathbf r)\mapsto\cX_k$, $1\le k\le N$.
A {\itbf Goursat problem} consists of prescribing the values $a_k({\bf x})=
a_{k0}({\bf x})$ on the subsets
\be\label{DGoursatData}
{\cal G}^{\bfeps}_k=\Big\{\sum_{i\in\cD_k}\mu_i\eps_i{\bf e}_i:\;\;
\mu_i\in{\mathbb Z},\;0\leq\mu_i\leq r_i/\eps_i\Big\}
\subset\Omega^{\bfeps}({\mathbf r}).
\ee
\end{dfn}
The dependent variables $a_k({\bf x})=a_k(x_1,\ldots,x_d)$ are thought of as
attached to the cubes of dimension $\#(\cD_k)$ adjacent to the point $\bf x$:
\be\label{cell}
c({\bf x};\cD_k):=
\Big\{{\bf x}+\sum_{i\in \cD_k}\mu_i\eps_i{\bf e}_i:\;\;0\le\mu_i\le 1\Big\}.
\ee
We use an abbreviation $\vec{f}_i(\vec{a})$ for the $N$--vector with the
components
\be
(\vec{f}_i(\vec{a}))_k=\left\{\begin{array}{cl}
        f_{(k,i)}(\vec{a})&\mbox{if } i\in{\cal E}_k,\\
        {\rm not\;\;defined}&\mbox{otherwise}.
\end{array}\right.
\ee
The {\itbf compatibility conditions} mentioned in the above definition
express the following requirement:
\be\label{comp}
\delta^{\eps_j}_{x_j}\delta^{\eps_i}_{x_i}a_k=
\delta^{\eps_i}_{x_i}\delta^{\eps_j}_{x_j}a_k
\ee
for any choice of $i\neq j$ from the respective ${\cal E}_k$. They read:
\begin{enumerate}
\item
The functions $f_{(k,i)}$ depend only on those $a_\ell$ for which
\be\label{comp1}
        {\cal E}_k\setminus\{i\}\subset{\cal E}_\ell\;.
\ee
\item For any $k=1,\ldots,N$ and any pair $i,j\in\cE_k$, the identity
\be\label{comp2}
\eps_if_{(k,i)}(\vec{a})+\eps_jf_{(k,j)}(\vec{a}+\eps_i\vec{f}_i(\vec{a}))
=\eps_jf_{(k,j)}(\vec{a})+\eps_if_{(k,i)}(\vec{a}+\eps_j\vec{f}_j(\vec{a}))
\ee
holds identically in $\vec{a}\in\cX_1\times\ldots\times\cX_N$.
\end{enumerate}
Indeed, in order for Eq. (\ref{comp}) to make sense, it is necessaary that
its both sides are well defined. In order for the left--hand side, say, to be
defined, the function $f_{(k,i)}$ is allowed to depend only on those
$a_{\ell}$ for which $j\in\cE_{\ell}$.
This has to hold for all $j\in\cE_k$, $j\neq i$, so we
come to (\ref{comp1}). Then (\ref{comp2}) is nothing but the in--length
translation of (\ref{comp}).
\begin{prp}
A Goursat problem for a compatible discrete hyperbolic system admits a unique
solution on $\Omega^{\bfeps}({\bf r})$.
\end{prp}

When a discrete hyperbolic system is considered as a discretization of a
continuous system supplied by B\"acklund transformations, it is naturally
supposed that the independent variables are divided into $(x_1,\ldots,x_n)$
discretized with the step $\eps_1=\ldots=\eps_n=\eps$, while the rest ones
$(x_{n+1},\ldots,x_d)$ are intrinsically discrete with
$\eps_{n+1}=\ldots=\eps_d=1$. We will write in this case
$\Omega^\eps(\mathbf r)$ for $\Omega^{\bfeps}(\mathbf r)$, so that in
the continuous limit $\eps\to 0$ the domain of independent variables becomes
\[
\Omega^0(\mathbf r)=[0,r_1]\times\ldots[0,r_n]
\times\{0,1,\ldots,r_{n+1}\}\times\ldots\times\{0,1,\ldots,r_d\} .
\]
The functions $f_{(k,i)}=f_{(k,i)}^\eps$ are
supposed to depend continuously on $\eps\in[0,\eps_0]$. In the limit
$\eps\to 0$ the first $\le n$ of equations (\ref{GDDH}) will turn into
differential ones:
\be\label{lim cont}
\partial_{x_i}a_k=f^0_{(k,i)}(\vec{a}), \quad i\in\cE_k,\quad 1\le i\le n,
\ee
while the rest ones will remain difference:
\be\label{lim discr}
\delta_{x_i}a_k=f^0_{(k,i)}(\vec{a}), \quad i\in\cE_k,\quad n+1\le i\le d.
\ee
The sets (\ref{DGoursatData}) on which the Goursat data are prescribed
will turn into
\be\label{CGoursatData}
{\cal G}_k^0=\Big\{\sum_{i\in\cD_k,\,i\le n}\mu_i{\bf e}_i+
\sum_{i\in\cD_k,\,i>n}\nu_i{\bf e}_i:\;\;
\mu_i\in[0,r_i],\,\nu_i\in\{0,1,\ldots,r_i\}\Big\}.
\ee

\begin{thm}
Let there be given an $\eps$--family of Goursat problems for compatible
discrete hyperbolic systems (\ref{GDDH}) on $\Omega^\eps(\mathbf r)$; denote
their solutions by $\vec{a}^{\,\eps}$. Suppose that
\be\label{conv f}
f^\eps_{(k,i)}(\vec{a})=f^0_{(k,i)}(\vec{a})+\cO(\eps)
\ee
uniformly on any compact subset of $\cX_1\times\ldots\times\cX_N$, and that
\be\label{conv a}
a_{k0}^\eps({\bf x})=a_{k0}^0({\bf x})+\cO(\eps)
\ee
uniformly on $\cG_k^\eps$. Then there exist $\bar{r}_i\in(0,r_i]$ for
$1\le i\le n$ and Lipschitz--continuous functions $\vec{a}^{\,0}$ on
$\Omega^0(\bar{r}_1,\ldots,\bar{r}_n,r_{n+1},\ldots,r_d)$ such that
\be\label{conv sol}
\vec{a}^{\,\eps}={\vec{a}}^{\,0}+\cO(\eps),
\ee
and $\vec{a}^0$ constitute the unique solution of the continuos Goursat
problem for the system (\ref{lim cont}), (\ref{lim discr}) with the
Goursat data on (\ref{CGoursatData}). If, in addition,
\begin{itemize}
\item the convergence (\ref{conv f}) is locally uniform in $C^{K+1}$,
\item the difference quotients of order $\le K+1$ of the discrete Goursat
data are bounded independently of $\eps$:
\[
\left|\Big(\prod_{i=1}^{K+1}\delta_{x_{j_i}}^\eps\Big)
a_{k0}^{\eps}\right|\le M\quad on\quad \cG_k^\eps,
\]
where all $j_i$ in the product have to belong to $\cD_k\cap\{1,\ldots,n\}$,
\item and all limit functions $a_{k0}^0$ belong to $C^K$,
\end{itemize}
then the convergence (\ref{conv sol}) is in the sense of $C^K$, i.e.
\[
\sup\left|\Big(\prod_{i=1}^K\delta^{\eps}_{x_{j_i}}\Big)a^\eps_k({\bf x})-
\Big(\prod_{i=1}^K\partial_{x_{j_i}}\Big)a_k^0({\bf x})\right|\to 0\quad
\mbox{ as }\quad\eps\rightarrow 0,
\]
where all $j_i$ in the products are in $\{1,\ldots,n\}$,
and the supremum is taken over
$\Omega^\eps(\bar{r}_1-K\eps,\ldots,\bar{r}_n-K\eps,r_{n+1},\ldots,r_d)$.
\end{thm}
{\bf Proof} of this theorem is completely analogous to the proofs of the
two-- and three--dimensional results given above. \BBox

\section{Appendix}
\label{Sect App}

{\bf Proof of Lemma \ref{Gronwall}.} Set $p=p(x)=\sum_{j=1}^d(x_j/\eps_j)$.
Suppose that the
statement is valid for all $x$ with $p(x)<p_0$, and consider some $x$
with $p(x)=p_0$. Notice that $p(x-\eps_i\bfe_i)=p(x)-1$. Therefore,
\[
\Delta(x)\le \max\Big(\Delta(0),\frac{\kappa}{\cal K}\Big)(1+\eps_i\cK)
\exp\Big(2{\cK}\sum_{j=1}^d x_j-2\cK\eps_i\Big)+\eps_i\kappa.
\]
It remains to estimate the last term on the right--hand side by
\beq
\eps_i\kappa & \le &
\frac{\kappa}{\cal K}\Big(e^{2\cK\eps_i}-(1+\eps_i\cK)\Big)\\
& \le & \max\Big(\Delta(0),\frac{\kappa}{\cal K}\Big)
\Big(e^{2\cK\eps_i}-(1+\eps_i\cK)\Big)
\exp\Big(2{\cK}\sum_{j=1}^d x_j-2\cK\eps_i\Big).
\eeq
Induction with respect to $p_0$ proves the Lemma. \BBox
\vfree

{\bf Proof of Lemma \ref{Gron}.} We prove first the case $Q=0$.
Since $\Delta$ is continuous and
hence bounded, there exists some $M\geq 0$ such that
\be\label{Gron aux1}
\Delta(x_1,\ldots,x_d)\leq M\exp\Big(2dL\sum_{i=1}^d x_i\Big)
\quad{\rm on}\quad \Omega({\mathbf r}).
\ee
We show that if this inequality holds for some $M>0$, then it holds also with
$M$ replaced by $M/2$. Indeed, (\ref{Gron aux1}) yields:
\beq
\Delta(x_1,\ldots,x_d) & \leq & LM\exp\Big(2dL(x_1+\ldots+x_d)\Big)
\sum_{j=1}^d\int_0^{x_j}e^{2dL(\xi_j-x_j)}d\xi_j  \\
 & = &
LM\exp\Big(2dL(x_1+\ldots+x_d)\Big)\sum_{j=1}^d \frac{1-\exp(-2dLx_j)}
{2dL}\\
 & \leq & (M/2)\exp\Big(2dL(x_1+\ldots+x_d)\Big).
\eeq
We conclude that (\ref{Gron aux1}) holds with any $M>0$, and so
$\Delta\leq 0$. Next, consider the case $Q>0$. Introduce the auxiliary
function
\[
q(x_1,\ldots,x_d):=2Q\exp\Big(2dL\sum_{i=1}^d x_i\Big).
\]
Exactly as above, we show:
\beq
L\sum_{j=1}^d\int_0^{x_j}
q(x_1,\ldots,x_{j-1},\xi_j,x_{j+1},\ldots,x_d)d\xi_j & \leq &
\frac{1}{2}q(x_1,\ldots,x_d) \\
 & \leq & q(x_1,\ldots,x_d)-Q.
\eeq
Therefore, the function $\Delta-q$ satisfies (\ref{GronwallX}) with $Q=0$,
so that $\Delta-q\leq 0$. \BBox
\vfree

{\bf Proof of Lemma \ref{composition}.}
Expand $f$ with respect to $a$ and $b$ around $a_0=a(x_0,y_0)$ and
$b_0=b(x_0,y_0)$:
\be
\label{conc4}
f(a,b)=\sum_{k+\ell\le K}{\cal L}_{k\ell}[(a-a_0)^{\otimes k},
(b-b_0)^{\otimes\ell}]+\rho,
\ee
where ${\cal L}_{k\ell}:=(k!\ell!)^{-1}D_a^kD_b^\ell f(a_0,b_0)$ is
a $(k+\ell)$-linear symmetric map.
By Taylor's theorem,
\beqn\label{roest}
|\rho| & \le & \sum_{k+\ell=K+1}\nn(D^k_aD^\ell_bf)[a,b]\nn_0
\cdot(|a-a_0|+|b-b_0|)^{K+1} \nonumber\\
 & \le &
\sum_{k+\ell=K+1}\nn(D^k_aD^\ell_bf)[a,b]\nn_0\cdot(L(|x-x_0|+|y-y_0|))^{K+1},
\eeqn
where $L$ is a Lipschitz constant of $a$ and $b$. To apply
$(\delta_x^\eps)^m(\delta^\eps_y)^n$ to $\rho$ means to evaluate
a weighted sum (with $\eps$--independent weights) of $\rho(x,y)$ at
lattice sites $(x,y)$ no more than $K$ steps away from $(x_0,y_0)$,
and then to divide by $\eps^K$. Using the estimate (\ref{roest}),
we find with some $C_K>0$:
\beq
|(\delta_x^\eps)^m(\delta^\eps_y)^n\rho(x_0,y_0)|
    &\le& \eps^{-K}C_K\sum_{|(x,y)-(x_0,y_0)|\le K\eps}|\rho(x,y)|\\
    &\le& \eps^{-K}C_KK^2 \sum_{k+\ell=K+1}\nn(D^k_aD^\ell_bf)[a,b]\nn_0
        \cdot (LK\eps)^{K+1}\\
    &\le& \eps\cdot C_K K^{K+3} L^{K+1}
    \sum_{k+\ell=K+1}\nn(D^k_aD^\ell_bf)[a,b]\nn_0:=B\eps.
\eeq
This is the $B$-term in Lemma \ref{composition}.

Next, we apply $(\delta_x^\eps)^m(\delta^\eps_y)^n$ to the sum
on the right--hand side of (\ref{conc4}). We have to estimate
expressions like $(\delta^\eps_x)^m(\delta^\eps_y)^n{\cal L}[c_1,\ldots,c_k]$,
where ${\cal L}[c_1,\ldots,c_k]$ denotes some $k$-linear map from $\cX^k$ to
$\cX$. Let ${\cal L}^*(C)$ with $C=\{c_1,\ldots,c_P\}$
symbolically stand for any linear combination of such expressions
\be
\sum_{j=1}^J \lambda_j\cdot {\cal L}[c_{j,1},\ldots,c_{j,k}]
\ee
with arguments $c_{j,i}$ arbitrarily chosen from $C$.

Assume that $c_i$ are actually functions on $\Omega^\eps(r)$,
so that ${\cal L}[c_1,\ldots,c_k]$ is $(x,y)$-dependent.
We start with the observation that
\be\label{conc1}
(\delta_x^\eps){\cal L}[c_1,\ldots,c_k]=
        \sum_{\kappa=1}^k{\cal L}[\tilde{c}_1,\ldots,\tilde{c}_{\kappa-1},
              \delta_x^\eps c_\kappa,\tilde{c}_{\kappa+1},\ldots,\tilde{c}_k]
\ee
where each $\tilde{c}_i$ stands either for the function $c_i(x,y)$
itself or for the shifted one $c_i(x+\eps,y)$. For instance,
one possible choice is $\tilde{c}_i=c_i(x,y)$ for $i<\kappa$ and
$\tilde{c}_i=c_i(x+\eps,y)$ for $i>\kappa$ in the $\kappa$th term.
The equality (\ref{conc1}) parallels the Leibnitz rule and follows
immediately from the multilinearity of ${\cal L}$.
By induction, we conclude that
\beqn
\label{conc2}
(\delta^\eps_x)^m{\cal L}[c_1,\ldots,c_k]&=&\sum_{\kappa=1}^k{\cal L}
   [\tilde{c}_1,\ldots,\tilde{c}_{\kappa-1},
   (\delta_x^\eps)^mc_\kappa,\tilde{c}_{\kappa+1},\ldots,\tilde{c}_k]\nonumber\\
        &&+{\cal L}^*\left(\Big\{(\delta_x^\eps)^i\tilde{c}_\kappa\Big\}_
                {{i<m}\atop{\kappa=1,\ldots,k}}\right)
\eeqn
where the tilde now denotes a shift in $x$ by no more than $m\eps$.
A further induction allows us to extend (\ref{conc2}) to the case of mixed
partial difference quotients:
\beqn
\label{conc3}
(\delta^\eps_x)^m(\delta^\eps_y)^n
        {\cal L}[c_1,\ldots,c_k]&=&\sum_{\kappa=1}^k{\cal L}
        [\tilde{c}_1,\ldots,\tilde{c}_{\kappa-1},
              (\delta_x^\eps)^m(\delta^\eps_y)^nc_\kappa,
              \tilde{c}_{\kappa+1},\ldots,\tilde{c}_k]\nonumber\\
  &&+{\cal L}^*\left(\Big\{(\delta_x^\eps)^i(\delta^\eps_y)^j\tilde{c_\kappa}
  \Big\}_{{i+j<K}\atop{\kappa=1,\ldots,k}}\right).
\eeqn
The tilde now stands for a possible shift in $x$ and $y$ by no more than
$m\eps$ and $n\eps$, respectively.

Applying (\ref{conc3}) to (\ref{conc4}) (one has to set
$c_1=\cdots=c_k=a-a_0$ and $c_{k+1}=\cdots=c_{k+\ell}=b-b_0$), we find:
\beqn\label{conc5}
k{\cal L}_{k\ell}[(\delta_x^\eps)^m(\delta^\eps_y)^na,
        (\tilde{a}-a_0)^{\otimes(k-1)},(\tilde{b}-b_0)^{\otimes\ell}]+
    \nonumber\\
\ell {\cal L}_{k\ell}[(\delta_x^\eps)^m(\delta^\eps_y)^nb,
     (\tilde{a}-a_0)^{\otimes k},(\tilde{b}-b_0)^{\otimes(\ell-1)}]+
     \nonumber\\
  +{\cal L}^*_{k\ell}\left(\{(\delta_x^\eps)^i(\delta^\eps_y)^j\tilde{a},
     (\delta_x^\eps)^i(\delta^\eps_y)^j\tilde{b}\}_{i+j<m+n}\right).
\eeqn
Taking into account that $|a-a_0|, |b-b_0|\le KL\eps$, we can estimate
the first two lines of (\ref{conc5}) as
\[
A\cdot\left(|(\delta^\eps_x)^m(\delta^\eps_y)^na(x,y)|+
|(\delta^\eps_x)^m(\delta^\eps_y)^na(x,y)|\right),
\]
where
\[
A=\sum_{0<k+\ell\le K}\frac{k+\ell}{k!\ell!}(KL\eps)^{k+\ell-1}
    \nn (D^k_aD^\ell_bf)[a,b]\nn_0.
\]
This is the $A$-term in Lemma \ref{composition}. Finally, the last line
of (\ref{conc5}) gives rise to the polynomial expression $P$ of
Lemma \ref{composition}:
\[
P(s,t)=C_K\sum_{0<k+\ell\le K}\nn (D^k_aD^\ell_bf)[a,b]\nn_0s^kt^\ell
\]
after estimating $|(\delta_x^\eps)^i(\delta^\eps_y)^j\tilde{a}(x_0,y_0)|$ by
$\nn a\nn_{K-1}$ and similarly for $b$.
\BBox

{\bf Remark.}
An immediate consequence of the the formula (\ref{conc3}) is a
product rule for discrete $C^k$-norms: If a product is defined between
elements of $\cX$, we may choose ${\cal L}[c_1,c_2]:=c_1\cdot c_2$
and so conclude
\be
\nn c_1\cdot c_2\nn_K\le C_K\nn c_1\nn_K\cdot\nn c_2\nn_K,
\ee
with $C_K>0$ depending only on $K$ and properties of the product.

\end{document}